\documentclass[a4paper]{amsart}
\usepackage{amsfonts}
\usepackage{amssymb}
\usepackage{amsmath}  
\usepackage{verbatim}

\input xy
\xyoption{all}

\newtheorem{theorem}{Theorem}[section]
\newtheorem{lemma}[theorem]{Lemma}

\newtheorem{proposition}[theorem]{Proposition}
\newtheorem{corollary}[theorem]{Corollary}
\newtheorem{definition}[theorem]{Definition}
\newtheorem{remark}[theorem]{Remark}
\newtheorem{notation}[theorem]{Notation}

\newcommand{\Ker}{\mbox{\rm Ker\,}}
\newcommand{\Coker}{\mbox{\rm Coker\,}}
\newcommand{\Ima}{\mbox{\rm Im\,}}
\newcommand{\End}[1]{\mbox{\rm End}_{#1}}
\newcommand{\Hom}[1]{\mbox{{\rm Hom}}_{#1}}
\newcommand{\Ext}[2]{\mbox{\rm Ext}^{#1}_{#2}}

\newcommand{\dimv}{\underline{\dim}\,}

\newcommand{\per}{\mbox{{\rm per }}}
\newcommand{\pr}[1]{\mbox{{\rm pr}}_{#1}}
\newcommand{\add}{\mbox{{\rm add \!}}}
\newcommand{\MOD}{\mbox{{\rm mod }}}
\newcommand{\Mod}{\mbox{{\rm Mod }}}
\newcommand{\proj}{\mbox{{\rm proj }}}
\newcommand{\perf}{\mbox{{\rm per \!}}}
\newcommand{\Ho}{\mbox{{\rm H}}}
\newcommand{\ind}[1]{\mbox{{\rm ind}}_{#1}}

\newcommand{\Gr}[1]{\mbox{{\rm Gr}}_{#1}}

\newcommand{\K}{\mbox{{\rm K}}}
\newcommand{\colim}{\mbox{{\rm colim \!}}}
\newcommand{\Lim}{\mbox{{\rm lim \!}}}

\newcommand{\demo}[1]{\textsc{Proof} #1 \hfill $\Box$ \bigskip}

\newcommand{\cC}{\mathcal{C}}
\newcommand{\cD}{\mathcal{D}}

\newcommand{\cF}{\mathcal{F}}
\newcommand{\cR}{\mathcal{R}}
\newcommand{\cT}{\mathcal{T}}
\newcommand{\cU}{\mathcal{U}}
\newcommand{\bQ}{\mathbb{Q}}

\begin{document}

\title[Cluster characters for any cluster category]{Cluster characters for cluster categories with infinite-dimensional morphism spaces}
\author{Pierre-Guy Plamondon}
\address{Universit\'e Paris Diderot -- Paris 7\\
   Institut de Math\'ematiques de Jussieu, UMR 7586 du CNRS \\
   Case 7012\\
   B\^atiment Chevaleret\\
   75205 Paris Cedex 13\\
   France }
\email{plamondon@math.jussieu.fr}
\thanks{The author was financially supported by an NSERC scholarship.}

\begin{abstract}
We prove the existence of cluster characters for $\Hom{}$-infinite cluster categories.  For this purpose, we introduce a suitable mutation-invariant subcategory of the cluster category.  We sketch how to apply our results in order to categorify any skew-symmetric cluster algebra.  More applications and a comparison to Derksen-Weyman-Zelevinsky's results will be given in a future paper.
\end{abstract}

\maketitle

\tableofcontents

\section{Introduction}
In their series of papers \cite{FZ02}, \cite{FZ03}, \cite{BFZ05}  and \cite{FZ07} published between 2002 and 2007, S.~Fomin and A.~Zelevinsky, together with A.~Berenstein for the third paper, introduced and developped the theory of cluster algebras.  They were motivated by the search for a combinatorial setting for total positivity and canonical bases.  Cluster algebras are a class of commutative algebras endowed with a distinguished set of generators, the cluster variables.  The cluster variables are grouped into finite subsets, called clusters, and are defined recursively from initial variables by repeatedly applying an operation called mutation on the clusters.  Recent surveys of the subject include \cite{Zelevinsky07}, \cite{GLS07} and \cite{Keller09}.

Cluster categories were introduced by A.~Buan, R.~Marsh, M.~Reineke, I.~Reiten and G.~Todorov in \cite{BMRRT06}, and by P.~Caldero, F.~Chapoton and R.~Schiffler in \cite{CCS06} for the $A_n$ case, in order to give a categorical interpretation of mutation of cluster variables.  In \cite{CC06}, P.~Caldero and F.~Chapoton used the geometry of quiver Grassmannians to define a map which, as they showed, yields a bijection from the set of isomorphism classes of indecomposable objects of the cluster category of a Dynkin quiver to the set of cluster variables in the associated cluster algebra.  It was proved by P.~Caldero and B.~Keller in \cite{CK06} that, for cluster algebras associated with acyclic quivers, the Caldero-Chapoton map induces a bijection between the set of isomorphism classes of indecomposable rigid objects and the set of cluster variables.

Using the notion of quiver with potential as defined in \cite{DWZ08}, C.~Amiot generalized the definition of cluster category in \cite{Amiot08}.  In the case where the quiver with potential is Jacobi-finite, the cluster character of Y.~Palu introduced in \cite{Palu08} sends reachable indecomposable rigid objects of the (generalized) cluster category to cluster variables.

Another approach for categorification of cluster algebras is studied by C.~Geiss, B.~Leclerc and J.~Schr\"oer in \cite{GLS05}, \cite{GLS06}, \cite{GLS072} and \cite{GLS08} where the authors use the category of modules over preprojective algebras of acyclic type.

In both cases, the categories encountered enjoy the following properties:  (1) they are $\Hom{}$-finite, meaning that the spaces of morphisms between any two objects is finite-dimensional; and (2) they are $2$-Calabi--Yau in the sense that for any two objects $X$ and $Y$, there is a bifunctorial isomorphism
\begin{displaymath}
	\Ext{1}{}(X, Y) \cong D\Ext{1}{}(Y, X).
\end{displaymath}

In this paper, we study a version of Y.~Palu's cluster characters for $\Hom{}$-infinite cluster categories, that is, cluster categories with possibly infinite-dimensional morphism spaces.  This cluster character $L \mapsto X'_L$ is not defined for all objcts $L$ but only for those in a suitable subcategory $\cD$, which we introduce.  We show that $\cD$ is mutation-invariant (in a sense to be defined) and that, for any objects $X$ and $Y$ of $\cD$, there is a bifunctorial non-degenarate bilinear form
\begin{displaymath}
	\Hom{}(X,\Sigma Y) \times \Hom{}(Y, \Sigma X) \longrightarrow k
\end{displaymath}
(this can be thought of as an adapted version of the $2$-Calabi--Yau property).

The category $\cD$ is equivalent to a $k$-linear subcategory of a certain derived category (the analogue of C.~Amiot's fundamental domain $\cF$ in \cite{Amiot08}).  We show that this subcategory also enjoys a certain property of invariance under mutation, as was first formulated as a ``hope'' by K.~Nagao in \cite{Nag}.

The main feature of the definition of the subcategory $\cD$ is the requirement that for any object $X$ of $\cD$, there exists a triangle
\begin{displaymath}
	T_1^X \longrightarrow T_0^X \longrightarrow X \longrightarrow \Sigma T_1^X
\end{displaymath}
where $T_0^X$ and $T_1^X$ are direct sums of direct summands of a certain fixed rigid object $T$.  This allows a definition of the \emph{index} of $X$, as in \cite{DK08} and \cite{Palu08}.

The main result of this article, besides the definition and study of the subcategory $\cD$, is the proof of a multiplication formula analogous to that of \cite{Palu08}: if $X$ and $Y$ are two objetcts of $\cD$ such that the spaces $\Hom{}(X, \Sigma Y)$ and $\Hom{}(Y, \Sigma X)$ are one-dimensional, and if 
\begin{displaymath}
\xymatrix{X \ar[r] & E \ar[r] & Y \ar[r] & \Sigma X } \textrm{ and }
\xymatrix{Y \ar[r] & E' \ar[r] & X \ar[r] & \Sigma Y } 
\end{displaymath}
are two non-split triangles, then we have the equality
\begin{displaymath}
	X'_X X'_Y = X'_E + X'_{E'}.
\end{displaymath}

This cluster character is in particular defined for the cluster category of any non-degenerate quiver with potential in the sense of \cite{DWZ08} (be it Jacobi-finite or not), and thus gives a categorification of any skew-symmetric cluster algebra.  Applications to cluster algebras will be the subject of a subsequent paper by the author.

In a different setting, using decorated representations of quivers with potentials, a categorification of any skew-symmetric cluster algebra was obtained by H.~Derksen, J.~Weyman and A.~Zelevinsky in the papers \cite{DWZ08} and \cite{DWZ09}, and these results were used by the authors to prove almost all of the conjectures formulated in \cite{FZ07}.

The article is organized as follows.

In Section \ref{sect::clustercat}, the main results concerning cluster categories of a quiver with potential and mutation are recalled.  In particular, we include the interpretation of mutation as derived equivalence, after \cite{KY09}.  The subcategory $\pr{\cC}\Gamma$, needed to define the subcategory $\cD$, is introduced and studied from Subsection \ref{sect::prGamma} up to the next section.  In Subsection \ref{sect::nagao}, we prove a result on the mutation of objects in the derived category, confirming K.~Nagao's hope in \cite{Nag}.  With hindsight, a precursor of this result is \cite[Corollary 5.7]{IR07}.

Section \ref{sect::clusterchar} is devoted to the definition of the cluster character $X'_?$.  After some preliminary results, it is introduced in Subsection \ref{sect::deficlusterchar} together with the subcategory $\cD$.  The multiplication formula is then proved in Subsection \ref{sect::multiplication}.

Finally, a link with skew-symmetric cluster algebras is given in Section \ref{sect::application}.

Throughout the paper, the symbol $k$ will denote an algebraically-closed field.  When working with any triangulated category, we will use the symbol $\Sigma$ to denote its suspension functor.  An object $X$ of any such category is \emph{rigid} if the space $\Hom{\cC}(X, \Sigma X)$ vanishes.

\subsection*{Acknowledgements}

This work is part of my PhD thesis, supervised by Professor Bernhard Keller.  I would like to express here my gratitude for his patience and his enthusiasm in sharing his mathematical knowledge.  I would like to thank Yann Palu and Dong Yang for their comments on an earlier version of this paper, and for precious conversations on categorification of cluster algebras.  Finally, I thank Kentaro Nagao for his question in \cite{Nag} which lead to Theorem \ref{theo::mutation}.

\section{Cluster category}\label{sect::clustercat}

In this section, after a brief reminder on quivers with potentials, the cluster category of a quiver with potential is defined after \cite{Amiot08}.  Mutation in the cluster category is then recalled.  Finally, we construct a subcategory on which a version of the cluster character of \cite{Palu08} will be defined in Section \ref{sect::clusterchar}.

\subsection{Skew-symmetric cluster algebras}
We briefly review the definition of (skew-symmetric) cluster algebras (the original definition appeared in \cite{FZ02} using mutation of matrices; the use of quivers was described, for example, in \cite[Definition 7.3]{FZ03} in a slightly different way than the one used here, and in \cite[Section 1.1]{GSV03}).  This material will be used in Section \ref{sect::application}.

A \emph{quiver} is a quadruple $Q = (Q_0, Q_1, s, t)$ consisting of a set $Q_0$ of vertices, a set $Q_1$ of arrows, and two maps $s,t : Q_1 \longrightarrow Q_0$ which send each arrow to its source or target.  A quiver is \emph{finite} if it has finitely many vertices and arrows.

Let $Q$ be a finite quiver without oriented cycles of length at most $2$.  We will denote the vertices of $Q$ by the numbers $1,2,\ldots,n$.  Let $i$ be a vertex of $Q$.  One defines the \emph{mutation of $Q$ at $i$} to be the quiver $\mu_i(Q)$ obtained from $Q$ in three steps :
\begin{enumerate}
	\item for each subquiver of the form $\xymatrix{ j\ar[r]^{a} & i\ar[r]^{b} & \ell}$, add an arrow $[ba]$ from $j$ to $\ell$;
  \item for each arrow $a$ such that $s(a) = i$ or $t(a) = i$, delete $a$ and add an arrow $a^*$ from $t(a)$ to $s(a)$ (that is, in the opposite direction);
  \item delete the arrows of a maximal set of pairwise disjoint oriented cycles of length 2 (which may have appeared in the first step).
\end{enumerate}

A \emph{seed} is a pair $(Q,\bold{u})$, where $Q$ is a finite quiver without oriented cycles of length at most $2$, and $\bold{u} = (u_1, u_2, \ldots, u_n)$ is an (ordered) free generating set of $\mathbb{Q}(x_1, x_2, \ldots, x_n)$.  Recall that $n$ is the number of vertices of $Q$.

If $i$ is a vertex of $Q$, the \emph{mutation} of the seed $(Q,\bold{u})$ is a new seed $(Q', \bold{u}') = (u'_1, u'_2, \ldots, u'_n))$, where
\begin{itemize}
	\item $Q'$ is the mutated quiver $\mu_i(Q)$;
	\item $u'_j = u_j$ whenever $j\neq i$;
	\item $u'_i$ is given by the equality
	      \begin{displaymath}
	        u'_iu_i = \prod_{a\in Q_1, t(a) = i}x_{s(a)} + \prod_{b\in Q_1, s(b) = i}x_{t(b)}. 
        \end{displaymath}
\end{itemize}

\begin{definition}
  Let $Q$ be a finite quiver without oriented cycles of length at most $2$.  Define the \emph{initial seed} as the seed $(Q, \bold{x} = (x_1, x_2, \ldots, x_n))$.
  \begin{itemize}
	  \item A \emph{cluster} is any set $\bold{u}$ appearing in a seed $(R, \bold{u})$ obtained from the initial seed by a finite sequence of mutation.
	  \item A \emph{cluster variable} is any element of a cluster.
	  \item The \emph{cluster algebra} associated with $Q$ is the $\mathbb{Q}$-subalgebra of the field of rational functions $\mathbb{Q}(x_1, x_2, \ldots, x_n)$ generated by the set of all cluster variables.
  \end{itemize}
\end{definition}

\subsection{Quivers with potentials and Jacobi-finiteness}

We denote by $kQ$ the \emph{path algebra} of a finite quiver $Q$ (that is, the space of formal linear combinations of paths of $Q$ endowed with the obvious multiplication) and by $\widehat{kQ}$ the \emph{complete path algebra of Q} (in which infinite linear combinations of paths are allowed).  The latter is a topological algebra for the $\mathfrak{m}$-adic topology, where $\mathfrak{m}$ is the ideal generated by the arrows of $Q$.  A basic system of open neighborhoods of zero is given by the powers of $\mathfrak{m}$.

Following \cite{DWZ08}, we define a \emph{quiver with potential} as a pair $(Q,W)$, where $Q$ is a finite quiver and $W$ is an element of the space $ Pot(Q) = \widehat{kQ}/U$, where $U$ is the closure of the commutator $[\widehat{kQ},\widehat{kQ}]$.  In other words, $W$ is a (possibly infinite) linear combination of oriented cycles in $Q$, considered up to cyclic equivalence of cycles (see \cite[Definition 3.2]{DWZ08}).  The element $W$ is called a \emph{potential} on $Q$.  In this paper, the terms of any potential will always be cycles of lenght at least $2$.

For any arrow $a$ of $Q$, define the \emph{cyclic derivative} of $a$ as the continuous linear map $\partial_a$ from the space of potentials to $\widehat{kQ}$ acting as follows on (equivalence classes of) oriented cycles :
\begin{displaymath}
	\partial_a(b_r\ldots b_2b_1) = \sum_{b_i = a}b_{i-1}b_{i-2}\ldots b_1b_rb_{r-1}\ldots b_{i+1}.
\end{displaymath}

The \emph{Jacobian algebra} $J(Q,W)$ of a quiver with potential $(Q,W)$ is the quotient of $\widehat{kQ}$ by the closure of the two-sided ideal generated by the $\partial_a W$, where $a$ ranges over all arrows of $Q$.  In case $J(Q,W)$ is finite-dimensional, $(Q,W)$ is \emph{Jacobi-finite}.



\subsection{Mutation of quivers with potentials}\label{subs::mutationQP}

Before defining mutation of quivers with potentials, we must say a word on the process of reduction.  Two quivers with potentials $(Q,W)$ and $(Q',W')$ are \emph{right-equivalent} if $Q_0 = Q'_0$ and there exists an $R$-algebra isomorphism $\varphi : \widehat{kQ} \longrightarrow \widehat{kQ'}$, where $R = \bigoplus_{i \in Q_0}ke_i$, such that $\varphi(W)$ equals $W'$ in $Pot(Q')$.  In that case, it is shown in \cite{DWZ08} that the Jacobian algebras of the two quivers with potential are isomorphic.

Let $(Q,W)$ be a quiver with potential.  We say that $(Q,W)$ is \emph{trivial} when $W$ is a combination of cycles of length at least $2$ and $J(Q,W)$ is isomorphic to $R$.  We say it is \emph{reduced} when $W$ has no terms which are cycles of length at most $2$.  

The direct sum of two quivers with potentials $(Q,W)$ and $(Q',W')$ such that $Q_0 = Q'_0$ is defined as being $(Q'', W + W')$, where $Q''$ is the quiver with the same set of vertices as $Q$ and $Q'$ and whose set of arrows is the union of those of $Q$ and $Q'$.

\begin{theorem}[\cite{DWZ08}, Theorem 4.6 and Proposition 4.5]
Any quiver with potential $(Q,W)$ is right equivalent to a direct sum of a reduced one $(Q_{red}, W_{red})$ and a trivial one $(Q_{triv}, W_{triv})$, both unique up to right-equivalence.  Moreover, $J(Q,W)$ and $J(Q_{red}, W_{red})$ are isomorphic.
\end{theorem}

The quiver with potential $(Q_{red}, W_{red})$ is the \emph{reduced part} of $(Q,W)$.

We can now define the mutation of a quiver with potential $(Q,W)$.  Let $i$ be a vertex of $Q$ not involved in an oriented cycle of length $\leq 2$. We can assume that $W$ is written as a series of oriented cycles whose source and target are not both $i$ (by replacing cycles by cyclically-equivalent ones if necessary).  The \emph{mutation at the vertex $i$} is the new quiver with potential $\mu_i(Q,W)$ obtained from $(Q,W)$ in the following way:
\begin{enumerate}
  \item for each subquiver of the form $\xymatrix{ j\ar[r]^{a} & i\ar[r]^{b} & \ell}$, add an arrow $[ba]$ from $j$ to $\ell$;
  \item for each arrow $a$ such that $s(a) = i$ or $t(a) = i$, delete $a$ and add an arrow $a^*$ from $t(a)$ to $s(a)$ (that is, in the opposite direction).  This gives a new quiver $\tilde{Q}$;
  \item replace $W$ by a potential $\tilde{W} = [W] + \sum_{a,b \in Q_1, t(a) = s(b) = i}a^*b^*[ba]$, where $[W]$ is obtained from $W$ by substuting $ba$ for $[ba]$ in its terms every time $t(a) = s(b) = i$.
\end{enumerate}
We use the notation $\tilde{\mu}_i(Q,W) = (\tilde{Q}, \tilde{W})$.  The mutation $\mu_i(Q,W)$ of $(Q,W)$ at $i$ is then defined to be the reduced part of $\tilde{\mu}_i(Q,W)$.

Remark that $\mu_i(Q,W)$ \emph{can} have oriented cycles of length $2$, even if $(Q,W)$ did not have any.  This forbids us to make iterated mutation at an arbitrary sequence of vertices.  

A sequence $(i_1, \ldots, i_r)$ of vertices is \emph{admissible} if all the mutations 
\begin{displaymath}
(Q,W), \ \mu_{i_1}(Q,W), \ \mu_{i_2}\mu_{i_1}(Q,W), \ \ldots, \ \mu_{i_r}\ldots\mu_{i_2}\mu_{i_1}(Q,W)
\end{displaymath}
are defined, that is, if $i_m$ is not involved in an oriented cycle of length $2$ in $\mu_{i_{m-1}}\ldots\mu_{i_1}(Q,W)$, for $2\leq m \leq r$, or in $(Q,W)$ for $m=1$.

A quiver with potential is \emph{non-degenerate} if any sequence of vetices is admissible.  Since we work over an algebraically closed field, the following existence result applies.

\begin{proposition}[\cite{DWZ08}, Corollary 7.4]
  Suppose that $Q$ is a finite quiver without oriented cycles of length at most $2$.  If the field $k$ is uncountable, then there exists a potential $W$ on $Q$ such that $(Q,W)$ is non-degenerate.
\end{proposition}

\subsection{Complete Ginzburg dg algebras}\label{subs::dgalg}

Let $(Q,W)$ be a quiver with potential.  Following Ginzburg in \cite{G06}, we construct a differential graded (dg) algebra $\Gamma = \Gamma_{Q,W}$ as follows.

First construct a new graded quiver $\overline{Q}$ from $Q$.  The vertices of $\overline{Q}$ are those of $Q$; its arrows are those of $Q$ (these have degree $0$), to which we add
\begin{itemize}
	\item for any arrow $a: i\longrightarrow j$ of $Q$, an arrow $a^*:j\longrightarrow i$ of degree $-1$;
	\item for any vertex $i$ of $Q$, a loop $t_i : i\longrightarrow i$ of degree $-2$.
\end{itemize}

Then, for any integer $i$, let 
\begin{displaymath}
	\Gamma^i = \prod_{\omega {\scriptsize \textrm{ path of degree }}i}k\omega.
\end{displaymath} 
This defines the graded $k$-algebra structure of $\Gamma$.  Its differential $d$ is defined from its action on the arrows of $\overline{Q}$.  We put
\begin{itemize}
	\item $d(a) = 0$, for each arrow $a$ of $Q$;
	\item $d(a^*) = \partial_a W$, for each arrow $a$ of $Q$;
	\item $d(t_i) = e_i\Big(\sum_{a\in Q_1}(aa^* - a^*a)\Big)e_i$, for each vertex $i$ of $Q$.
\end{itemize}

The differential graded algebra $\Gamma$ thus defined is the \emph{complete Ginzburg dg algebra} of $(Q,W)$.  It is linked to the Jacobian algebra of $(Q,W)$ as follows.

\begin{lemma}[\cite{KY09}, Lemma 2.8]
With the above notations, $J(Q,W)$ is isomorphic to $\Ho^0\Gamma$.
\end{lemma}

\subsection{Cluster category}\label{subs::clustercat}
Keep the notations of Section \ref{subs::dgalg}.

Denote by $\cD\Gamma$ the derived category of $\Gamma$ (see \cite{K94} or \cite{KY09} for background material on the derived category of a dg algebra).  Consider $\Gamma$ as an object of $\cD\Gamma$.  The \emph{perfect derived category} of $\Gamma$ is the smallest full triangulated subcategory of $\cD\Gamma$ containing $\Gamma$ and closed under taking direct summands.  It is denoted by $\perf \Gamma$.

Denote by $\cD_{fd}\Gamma$ the full subcategory of $\cD\Gamma$ whose objects are those of $\cD\Gamma$ with finite-dimensional total homology. This means that homology is zero except in finitely many degrees, where it is of finite dimension.  As shown in \cite[Theorem 2.17]{KY09}, the category $\cD_{fd}\Gamma$ is a triangulated subcategory of $\perf \Gamma$.

Moreover, we have the following \emph{relative $3$-Calabi--Yau property} of $\cD_{fd}\Gamma$ in $\cD\Gamma$.

\begin{theorem}[\cite{K08}, Lemma 4.1 and \cite{K09}, Theorem 6.3]\label{theo::3CY}
  For any objects $L$ of $\cD\Gamma$ and $M$ of $\cD_{fd}\Gamma$, there is a canonical isomorphism
  \begin{displaymath}
	  D\Hom{\cD\Gamma}(M, L) \longrightarrow \Hom{\cD\Gamma}(\Sigma^{-3}L, M)
  \end{displaymath}
functorial in both $M$ and $L$.  
\end{theorem}

Following \cite[Definition 3.5]{Amiot08} (and \cite[Section 4]{KY09} in the non Jacobi-finite case), we define the \emph{cluster category} of $(Q,W)$ as the idempotent completion of the triangulated quotient $(\perf\Gamma) / \cD_{fd}\Gamma$, and denote it by $\cC = \cC_{Q,W}$.

In case $(Q,W)$ is Jacobi-finite, $\cC_{Q,W}$ enjoys the following properties (\cite[Theorem 3.6]{Amiot08} and \cite[Proposition 2.1]{KR07}) :
\begin{itemize}
	\item it is $\Hom{}$-finite;
	\item it is $2$-Calabi--Yau;
	\item the object $\Gamma$ is \emph{cluster-tilting} in the sense that it is rigid and any object $X$ of $\cC$ such that $\Hom{\cC}(\Gamma, \Sigma X) = 0$ is in $\add \Gamma$;
	\item any object $X$ of $\cC$ admits an $(\add \Gamma)$-\emph{presentation}, that is, there exists a triangle $\xymatrix{ T_1^X \ar[r] & T_0^X \ar[r] & X \ar[r] & \Sigma T_1^X}$, with $T_1^X$ and $T_0^X$ in $\add \Gamma$.
\end{itemize}

As we shall see later, most of these properties do not hold when $(Q,W)$ is not Jacobi-finite.

\subsection{Mutation in $\cC$}
Keep the notations of the previous section.  Let $i$ be a vertex of $Q$ not involved in any oriented cycle of length $2$.  As seen in Section \ref{subs::mutationQP}, one can mutate $(Q,W)$ at the vertex $i$.

In the cluster category, this corresponds to changing a direct factor of $\Gamma$.  Let $\Gamma'$ be the complete Ginzburg dg algebra of $\tilde{\mu}_i(Q,W)$.  For any vertex $j$ of $Q$, let $\Gamma_j = e_j\Gamma$ and $\Gamma'_j = e_j\Gamma'$.

\begin{theorem}[\cite{KY09}, Theorem 3.2]\label{theo::KY} 
 \begin{enumerate}
   \item There is a triangle equivalence $F$ from $\cD(\Gamma')$ to $\cD(\Gamma)$ sending $\Gamma'_j$ to $\Gamma_j$ if $i \neq j$ and to the cone $\Gamma_i^*$ of the morphism 
     \begin{displaymath}
	     \Gamma_i \longrightarrow \bigoplus_{\alpha}\Gamma_{t(\alpha)}
     \end{displaymath}
     whose components are given by left multiplication by $\alpha$ if  $i=j$.  The functor $F$ restricts to triangle equivalences from $\perf \Gamma'$ to $\perf \Gamma$ and from $\cD_{fd} \Gamma'$ to $\cD_{fd} \Gamma$.
     
   \item Let $\Gamma_{red}$ and $\Gamma'_{red}$ be the complete Ginzburg dg algebra of the reduced part of $(Q,W)$ and $\tilde{\mu}_i(Q,W)$, respectively.  The functor $F$ induces a triangle equivalence $F_{red} : \cD(\Gamma'_{red}) \longrightarrow \cD(\Gamma_{red})$ which restricts to triangle equivalences from $\perf \Gamma_{red}'$ to $\perf \Gamma_{red}$ and from $\cD_{fd} \Gamma'_{red}$ to $\cD_{fd} \Gamma_{red}$.
 \end{enumerate}
\end{theorem}

The object $\Gamma_i^* \oplus \bigoplus_{j \neq i} \Gamma_j$ in $\cD\Gamma$ is the \emph{mutation of $\Gamma$ at the vertex $i$}, and we denote it by $\mu_i(\Gamma)$.

Note that in part (2) of the above theorem, $e_j\Gamma'_{red}$ is still sent to $e_j\Gamma_{red}$ if $i \neq j$ and to the cone $\Gamma_{red, i}^*$ of the morphism 
     \begin{displaymath}
	     e_i\Gamma_{red} \longrightarrow \bigoplus_{\alpha} e_{t(\alpha)}\Gamma_{red}
     \end{displaymath}
     whose components are given by left multiplication by $\alpha$ if $i=j$.
     
For instance, if  $(i_1, i_2, \ldots, i_r)$ is an admissible sequence of vertices, then we get a sequence of triangle equivalences
\begin{displaymath}
	\cD\Gamma^{(r)} \longrightarrow \ldots \longrightarrow \cD\Gamma^{(1)} \longrightarrow \cD\Gamma,
\end{displaymath}
where $\Gamma^{(j)}$ is the complete Ginzburg dg-algebra of $\mu_{i_j}\mu_{i_{j-1}}\ldots\mu_{i_1}(Q,W)$, for $j\in \{ 1,2,    \ldots, r\}$.  We denote the image of $\Gamma^{(r)}$ in $\cD\Gamma$ by $\mu_{i_r}\mu_{i_{r-1}}\ldots\mu_{i_1}(\Gamma)$.

We now remark some consequences of \ref{theo::KY} on the level of cluster categories.  First, there are induced triangle equivalence $\cC_{\tilde{\mu}_i(Q,W)}\longrightarrow \cC_{Q,W}$ and $\cC_{\mu_i(Q,W)} \longrightarrow \cC_{Q,W}$.

Moreover, as shown in Section 4 of \cite{KY09}, the cone of the morphism
\begin{displaymath}
	\bigoplus_{\beta: i\rightarrow j}\Gamma_j \longrightarrow \Gamma_i
\end{displaymath}
whose components are given by left multiplication by $\beta$ is isomorphic to $\Sigma \Gamma_i^*$ in $\cC$.  Hence we have triangles in $\cC$
\begin{displaymath}
	\Gamma_i \longrightarrow \bigoplus_{\alpha : i\rightarrow j}\Gamma_j \longrightarrow \Gamma_i^* \longrightarrow \Sigma \Gamma_i \ \ \textrm{ and } \ \ \Gamma_i^* \longrightarrow \bigoplus_{\alpha : j\rightarrow i}\Gamma_j \longrightarrow \Gamma_i \longrightarrow \Sigma \Gamma_i^*,
\end{displaymath}
and $\dim \Hom{\cC}(\Gamma_j, \Sigma\Gamma_i^*) = \delta_{i,j}$ (see \cite[Section 4]{KY09}).

If $(Q,W)$ is non-degenerate and reduced, then any sequence of vertices $i_1, \ldots, i_r$ yields a sequence of triangle equivalences
\begin{displaymath}
	\cC_{\mu_{i_r} \ldots \mu_{i_1}(Q,W)} \longrightarrow \ldots \longrightarrow \cC_{\mu_{i_1}(Q,W)} \longrightarrow \cC_{Q,W}
\end{displaymath}
sending $\Gamma_{\mu_{i_r} \ldots \mu_{i_1}(Q,W)}$ to $\mu_{i_r} \ldots \mu_{i_1}(\Gamma_{Q,W})$.

\subsection{The subcategory $\pr{\cC}\Gamma$}\label{sect::prGamma}
  Since in general the cluster category does not enjoy the properties listed in Section \ref{subs::clustercat}, we will need to restrict ourselves to a subcategory of it.
  
  Let $\cT$ be any triangulated category.  For any subcategory $\cT'$ of $\cT$, define $\ind{} \cT'$ as the set of isomorphism classes of indecomposable objects of $\cT$ contained in $\cT'$.  Denote by $\add \cT'$ the full subcategory of $\cT$ whose objects are all finite direct sums of direct summands of objects in $\cT'$.  The subcategory $\cT'$ is \emph{rigid} if, for any two objects $X$ and $Y$ of $\cT'$, $\Hom{\cT}(X, \Sigma Y) = 0$.
  
  Finally, define $\pr{\cT}\cT'$ as the full subcategory of $\cT$ whose objects are cones of morphisms in $\add \cT'$ (the letters ``$\pr{}$'' stand for \emph{presentation}, as all objects of $\pr{\cT}\cT'$ admit an ($\add \cT'$)-presentation).  In the notations of \cite[Section 1.3.9]{BBD82}, this subcategory is written as $(\add \cT') \ast (\add \Sigma\cT')$.
  
  As we shall now prove, the category $\pr{\cT}\cT'$ is invariant under ``mutation'' of $\cT'$.
  
  Recall that a category $\cT'$ is \emph{Krull--Schmidt} if any object can be written as a finite direct sum of objects whose endomorphism rings are local.  Note that in that case, we have $\cT' = \add \cT'$.
  
\begin{proposition}\label{prop::prmutation}
  Let $\cR$ and $\cR'$ be rigid Krull--Schmidt subcategories of a triangulated category $\cT$.  Suppose that there exist indecomposable objects $R$ of $\cR$ and $R^*$ of $\cR'$ such that $\ind{}\cR \setminus \{ R \} = \ind{}\cR' \setminus \{ R^* \}$.  Suppose, furthermore, that $\dim \Hom{\cT}(R, \Sigma R^*) = \dim \Hom{\cT}(R^*, \Sigma R) = 1$.  Let 
  \begin{displaymath}
	  R \longrightarrow E \longrightarrow R^* \longrightarrow \Sigma R \ \ \textrm{ and } \ \ R^* \longrightarrow E' \longrightarrow R \longrightarrow \Sigma R^*
  \end{displaymath}
be non-split triangles, and suppose that $E$ and $E'$ lie in $\cR \cap \cR'$.

Then $\pr{\cT}\cR = \pr{\cT}\cR'$.
\end{proposition}
\demo{  In view of the symmetry of the hypotheses, we only have to prove that any object of $\pr{\cT}\cR$ is an object of $\pr{\cT}\cR'$.

Let $X$ be an object of $\pr{\cT}\cR$.  Let $T_1 \longrightarrow T_0 \longrightarrow X \longrightarrow \Sigma T_1$ be a triangle, with $T_1$ and $T_0$ in $\cR$.

The category $\cR$ being Krull--Schmidt, one can write (in a unique way up to isomorphism) $T_0 = \overline{T}_0 \oplus R^m$, where $R$ is not a direct summand of $\overline{T}_0$.

The composition $\overline{T}_0 \oplus (E')^m \longrightarrow \overline{T}_0 \oplus R^m \longrightarrow X$ yields an octahedron

$\phantom{an octahedron}$\begin{xy} 0;<1pt,0pt>:<0pt,-1pt>:: 
(105,0) *+{\Sigma W} ="0",
(74,90) *+{\overline{T}_0 \oplus (E')^m} ="1",
(207,90) *+{X.} ="2",
(0,69) *+{(\Sigma R^*)^m} ="3",
(103,147) *+{\overline{T}_0\oplus R^m} ="4",
(133,69) *+{\Sigma T_1} ="5",
"0", {\ar|+"1"},
"2", {\ar"0"},
"3", {\ar"0"},
"0", {\ar@{.>}"5"},
"1", {\ar"2"},
"3", {\ar|+"1"},
"1", {\ar"4"},
"4", {\ar"2"},
"2", {\ar@{.>}"5"},
"4", {\ar"3"},
"5", {\ar@{.>}|+"3"},
"5", {\ar@{.>}|+"4"},
\end{xy} 

Now write $T_1 = \overline{T}_1 \oplus R^n$.  Then we have a triangle
\begin{displaymath}
	\xymatrix{(R^*)^m \ar[r] & W\ar[r] & \overline{T}_1 \oplus R^n \ar[r]^{\varepsilon} & (\Sigma R^*)^m}.
\end{displaymath}
Since $\Hom{\cT}(\overline{T}_1, \Sigma R^*) = 0$ and $\dim \Hom{\cT}(R, \Sigma R^*) = 1$, by a change of basis, we can write $\varepsilon$ in matrix form as
\begin{displaymath}
	\left(\begin{array}{c|c}
    I_rx & 0 \\
    \hline
    0 & 0
  \end{array}\right)
\end{displaymath}
where $x$ is a non-zero element of $\Hom{\cT}(R, \Sigma R^*)$.  Therefore $W$ is isomorphic to $E^r \oplus R^{n-r} \oplus (R^*)^{m-r} \oplus \overline{T}_1$.  

Now, we have a triangle $W \longrightarrow \overline{T}_0 \oplus (E')^m \longrightarrow X \longrightarrow \Sigma W$.  Compose $\Sigma^{-1}X \longrightarrow W$ with $W = E^r \oplus R^{n-r} \oplus (R^*)^{m-r} \oplus \overline{T}_1 \longrightarrow E^r \oplus E^{n-r} \oplus (R^*)^{m-r} \oplus \overline{T}_1$ (the second term is changed) to get an octahedron

$\phantom{an octahed}$\begin{xy} 0;<1pt,0pt>:<0pt,-1pt>:: 
(105,0) *+{V} ="0",
(74,90) *+{\Sigma^{-1} X} ="1",
(207,90) *+{E^r \oplus E^{n-r} \oplus (R^*)^{m-r} \oplus \overline{T}_1.} ="2",
(0,69) *+{\overline{T}_0 \oplus (E')^m} ="3",
(103,147) *+{E^r \oplus R^{n-r} \oplus (R^*)^{m-r} \oplus \overline{T}_1} ="4",
(133,69) *+{(R^*)^{n-r}} ="5",
"0", {\ar|+"1"},
"2", {\ar"0"},
"3", {\ar"0"},
"0", {\ar@{.>}"5"},
"1", {\ar"2"},
"3", {\ar|+"1"},
"1", {\ar"4"},
"4", {\ar"2"},
"2", {\ar@{.>}"5"},
"4", {\ar"3"},
"5", {\ar@{.>}|+"3"},
"5", {\ar@{.>}|+"4"},
\end{xy}

The morphism $(R^*)^{n-r} \longrightarrow \overline{T}_0 \oplus (E')^m$ is zero, so the triangle $\overline{T}_0 \oplus (E')^m \longrightarrow V \longrightarrow (R^*)^{n-r} \longrightarrow \Sigma(\overline{T}_0 \oplus (E')^m)$ splits, and $V$ is isomorphic to $(R^*)^{n-r} \oplus \overline{T}_0 \oplus (E')^m$.

Hence we have a triangle
\begin{displaymath}
	E^n \oplus (R^*)^{m-r} \oplus \overline{T}_1 \longrightarrow (R^*)^{n-r} \oplus \overline{T}_0 \oplus (E')^m \longrightarrow X \longrightarrow \Sigma(E^n \oplus (R^*)^{m-r} \oplus \overline{T}_1),
\end{displaymath} 
proving that $X$ belongs to $\pr{\cT}\cR'$.  This finishes the proof.
}

\begin{corollary}\label{coro::pr}
Let $\cC$ be the cluster category of a quiver with potential $(Q,W)$.  For any admissible sequence $(i_1, \ldots, i_r)$ of vertices of $Q$, the following equality holds :
  \begin{displaymath}
	  \pr{\cC} \Gamma = \pr{\cC} \Big(\mu_{i_r}\ldots\mu_{i_1}(\Gamma)\Big).
  \end{displaymath}
\end{corollary}
\demo{ We apply Proposition \ref{prop::prmutation} and use induction on $r$. That $\add\Gamma$ is a Krull--Schmidt category is shown in Corollary \ref{coro::prkrull} below.  We also need that $\Gamma$ is a rigid object of $\cC$ ; this follows from Proposition \ref{prop::fund} below.
}

\subsection{Properties of $\pr{\cC}\Gamma$}

Let $\cC$ be the cluster category of a quiver with potential $(Q,W)$. We will prove in this section that the subcategory $\pr{\cC}\Gamma$ enjoys versions of some of the properties listed in Section \ref{subs::clustercat}.

We denote by $\cD_{\leq 0}$ (and $\cD_{\geq 0}$ respectively) the full subcategory of $\cD\Gamma$ whose objects are those $X$ whose homology is concentrated in non-positive (and non-negative, respectively) degrees.  Recall that $\cD_{\leq 0}$ and $\cD_{\geq 0}$ form a \emph{t-structure}; in particular, $\Hom{\cD\Gamma}(\cD_{\leq 0}, \cD_{\geq 1})$ vanishes, and for each object $X$ of $\cD\Gamma$, there exists a unique (up to a unique triangle isomorphism) triangle
\begin{displaymath}
	\tau_{\leq 0}X \longrightarrow X \longrightarrow \tau_{\geq 1}X \longrightarrow \Sigma \tau_{\leq 0}X
\end{displaymath}
with $\tau_{\leq 0}X$ in $\cD_{\leq 0}$ and $\tau_{\geq 1}X$ in $\cD_{\geq 1}$.

\begin{lemma}\label{lemm::fund}
If $X$ and $Y$ lie in $\pr{\cD\Gamma}\Gamma$, then the quotient functor $\perf \Gamma \longrightarrow \cC$ induces an isomorphism
\begin{displaymath}
	\Hom{\cD\Gamma}(X,Y) \longrightarrow \Hom{\cC}(X,Y).
\end{displaymath}
\end{lemma}
\demo{ Let $X$ and $Y$ be as in the statement. In particular, $X$ and $Y$ lie in $\cD_{\leq 0}\Gamma$.

First suppose that a morphism $f : X\longrightarrow Y$ is sent to zero in $\cC$.  This means that $f$ factors as 
\begin{displaymath}
 \xymatrix{X\ar[r]^g & M \ar[r]^h & Y}, 
\end{displaymath} 
with $M$ in $\cD_{fd}\Gamma$.  Now, $X = \tau_{\leq 1}X$, so $g$ factors through $\tau_{\leq 1}M$, which is still in $\cD_{fd}$.  Using Theorem \ref{theo::3CY}, we have an isomorphism
\begin{displaymath}
	D\Hom{\cD\Gamma}(\tau_{\leq 1}M, Y) \longrightarrow \Hom{\cD\Gamma}(Y, \Sigma^3 \tau_{\leq 1}M).
\end{displaymath}

The right hand side of this equation is zero, since $\Hom{\cD\Gamma}(Y, \cD_{\leq -2}\Gamma)=0$.  Hence $f=0$.  This shows injectivity.

To prove surjectivity, consider a fraction
\begin{displaymath}
	\xymatrix{X\ar[r]^{f} & Y' & Y \ar[l]_{s}},
\end{displaymath}
where the cone of $s$ is an object $N$ of $\cD_{fd}\Gamma$.

The following diagram will be helpful.
\begin{displaymath}
	\xymatrix{ & Y\ar@{=}[r]\ar[d]_s & Y\ar[d]_t \\
	           X\ar[r]^f & Y'\ar[r]^g\ar[d] & Y''\ar[d] \\
	           \tau_{\leq 0}N\ar[r] & N\ar[d]\ar[r] & \tau_{\geq 1}N \ar[d]_h \\
	           & \Sigma Y \ar@{=}[r] & \Sigma Y
	}
\end{displaymath}

We have that $\Hom{\cD\Gamma}(\tau_{\leq 0}N, \Sigma Y)$ is isomorphic to $D\Hom{\cD\Gamma}(Y, \Sigma^2 \tau_{\leq 0}N)$ because of Theorem \ref{theo::3CY}, and this space is zero since $\Hom{\cD\Gamma}(Y, \cD_{\leq -2})$ vanishes.  Thus there exists a morphism $h:\tau_{\geq 1}N \rightarrow \Sigma Y$ such that the lower right square of the above diagram commute.  We embed $h$ in a triangle; this triangle is the rightmost column of the diagram.

We get a new fraction
\begin{displaymath}
	\xymatrix{X\ar[r]^{gf} & Y'' & Y\ar[l]_{t} }
\end{displaymath}
which is equal to the one we started with.  But since $X$ is in $\cD_{\leq 0}$ and $\tau_{\geq 1}N$ is in $\cD_{\geq 1}$, the space $\Hom{\cD\Gamma}(X, \tau_{\geq 1}N)$ vanishes.  Thus there exists a morphism $\ell:X\rightarrow Y$ such that $gf = t\ell$.  It is easily seen that the fraction is then the image of $\ell$ under the quotient functor. 

Thus the map is surjective.
}

\begin{proposition}\label{prop::fund}
The quotient functor $\perf \Gamma \longrightarrow \cC$ restricts to an equivalence of ($k$-linear) categories $\pr{\cD\Gamma}\Gamma \longrightarrow \pr{\cC}\Gamma$.
\end{proposition}
\demo{ It is a consequence of Lemma \ref{lemm::fund} that the functor is fully faithful.

It remains to be shown that it is dense.  Let $Z$ be an object of $\pr{\cC}\Gamma$, and let $T_1\longrightarrow T_0 \longrightarrow Z \longrightarrow \Sigma T_1$ be an $\add \Gamma$-presentation.  The functor being fully faithful, the morphism $T_1\longrightarrow T_0$ lifts in $\pr{\cD\Gamma}\Gamma$ to a morphism $P_1\longrightarrow P_0$, with $P_0$ and $P_1$ in $\add\Gamma$.  Its cone is clearly sent to $Z$ in $\cC$.  This finishes the proof of the equivalence.
} 

As in \cite{Amiot08}, we have the following characterization of $\pr{\cD\Gamma}\Gamma$, which we shall prove after Corollary \ref{coro::prkrull}.

\begin{lemma}\label{rema::fund}
 We have that that $\pr{\cD\Gamma}\Gamma = \cD_{\leq 0} \cap {}^{\perp}\cD_{\leq -2} \cap \perf \Gamma$.
\end{lemma}

\begin{corollary}\label{coro::prkrull}
The category $\pr{\cC}\Gamma$ is a Krull--Schmidt category.
\end{corollary}
\demo{ In view of Proposition \ref{prop::fund}, it suffices to prove that $\pr{\cD\Gamma}\Gamma$ is a Krull--Schmidt category.  It is shown in \cite[Lemma 2.17]{KY09} that the category $\per\Gamma$ is a Krull--Schmidt category.  Since $\pr{\cD\Gamma}\Gamma$ is a full subcategory of $\per\Gamma$, it is sufficient to prove that any direct summand of an object in $\pr{\cD\Gamma}\Gamma$ is also in $\pr{\cD\Gamma}\Gamma$. The equality $\pr{\cD\Gamma}\Gamma = \cD_{\leq 0} \cap {}^{\perp}\cD_{\leq -2} \cap \perf \Gamma$ of Lemma \ref{rema::fund} implies this property.  Note that it also follows from \cite[Proposition 2.1]{IY08}, whose proof does not depend on the $\Hom{}$-finiteness assumption.
}

In order to prove Lemma \ref{rema::fund}, we will need the following definition.

\begin{definition}\label{defi::minimal}
A dg $\Gamma$-module $M$ is \emph{minimal perfect} if its underlying graded module is of the form
\begin{displaymath}
	\bigoplus_{j=1}^{N} R_j,
\end{displaymath}
where each $R_j$ is a finite direct sum of shifted copies of direct summands of $\Gamma$, and if its differential is of the form $d_{int} + \delta$, where $d_{int}$ is the direct sum of the differential of the $R_j$, and $\delta$, as a degree $1$ map from $\bigoplus_{j=1}^{N} R_j$ to itself, is a strictly upper triangular matrix whose entries are in the ideal of $\Gamma$ generated by the arrows.
\end{definition}

\begin{lemma}
Let $M$ be a dg $\Gamma$-module such that $M$ is perfect in $\cD\Gamma$.  Then $M$ is quasi-isomorphic to a minimal perfect dg module.
\end{lemma}
\demo{ We will apply results of \cite{Bondarko09}.  Using the notation of \cite[Section 6.2]{Bondarko09}, $\perf\Gamma$ is equivalent to the category Tr$(C)$, where $C$ is the dg category whose objects are vertices of the quiver $Q$ and morphisms dg vector spaces are given by the paths of $Q$.  Thus any object of $\perf\Gamma$ is quasi-isomorphic to a dg module as in Definition \ref{defi::minimal}, where the entries of $\delta$ do not necessarily lie in the ideal generated by the arrows.

As a graded $\Gamma$-module, any such object can be written in the form $\Sigma^{i_1}\Gamma_{j_1} \oplus \ldots \oplus \Sigma^{i_r}\Gamma_{j_r}$, where each $j_{\ell}$ is a vertex of $Q$ and each $i_{\ell}$ is an integer.  Assume that $i_1\leq \ldots \leq i_r$.  The subcategory of objects wich can be written in this form, with $a \leq i_1 \leq i_r \leq b$, is denoted by $\cC^{[a,b]}$.  According to \cite[Lemma 5.2.1]{Bondarko09}, $\cC^{[a,b]}$ is closed under taking direct summands.

Let $X$ be an object of $\perf\Gamma$.  Then there are integers $a$ and $b$ such that $X$ lies in $\cC^{[a,b]}$.  We prove the Lemma by induction on $b-a$.

If $a=b$, then $\delta$ has to be zero, and $X$ is minimal perfect.

Suppose that all objects of $\cC^{[a,b]}$ are isomorphic to a minimal perfect dg module whenever $b-a$ is less or equal to some integer $n\geq 0$.

Let $X$ be an object of $\cC^{[a,b]}$, with $b-a = n+1$.  We can assume that $X$ is of the form $\Sigma^{i_1}\Gamma_{j_1} \oplus \ldots \oplus \Sigma^{i_r}\Gamma_{j_r}$ and that its differential is written in matrix form as
\begin{displaymath}
  \left( \begin{array}{cccc}
    d_{\Sigma^{i_1}\Gamma_{j_1}} & f_{12} & \ldots & f_{1r} \\
    0 & d_{\Sigma^{i_2}\Gamma_{j_2}} & \ldots & f_{2r} \\
    \vdots & \vdots & \ddots & \vdots \\
    0 & 0 & \ldots & d_{\Sigma^{i_r}\Gamma_{j_r}}
  \end{array} \right),
\end{displaymath}
where all the $f_{uv}$ are in the ideal generated by the arrows.  

Suppose that $i_q = i_{q+1} = \ldots = i_r$, but $i_{q-1} < i_q$.  Then $X$ is the cone of the morphism from $\Sigma^{i_q-1}\Gamma_{j_q} \oplus \ldots \oplus \Sigma^{i_r-1}\Gamma_{j_r}$ to the submodule $X'$ of $X$ whose underlying graded module is $\Sigma^{i_1}\Gamma_{j_1} \oplus \ldots \oplus \Sigma^{i_{q-1}}\Gamma_{j_{q-1}}$ given by the matrix
  
\begin{displaymath}
  \left( \begin{array}{cccc}
    \Sigma^{-1}f_{1,q} & \Sigma^{-1}f_{1,q+1} & \ldots & \Sigma^{-1}f_{1,r} \\
    \Sigma^{-1}f_{2,q} & \Sigma^{-1}f_{2,q+1} & \ldots & \Sigma^{-1}f_{2,r} \\
    \vdots & \vdots &  & \vdots \\
    \Sigma^{-1}f_{q-1,q} & \Sigma^{-1}f_{q-1,q+1} & \ldots & \Sigma^{-1}f_{q-1,r}
  \end{array} \right),
\end{displaymath}
whose entries are still elements of $\Gamma$.  Note that $X'$ lies in $\cC^{[a,b-1]}$.  By the induction hypothesis, $X'$ is quasi-isomorphic to a minimal perfect dg module.  Thus we can assume that $f_{ij}$ is in the ideal generated by the arrows, for $i = 1,2, \ldots q-1$ and $j= 1,2,\ldots, q-1$.

The rest of the proof is another induction, this time on the number of summands of $X$ of the form $\Sigma^m\Gamma_{\ell}$ (this number is $r-q+1$).  

If this number is $1$, then $X$ is the cone of a morphism given in matrix form by a column.  If this column contains no isomorphisms, then $X$ is minimal perfect.  Otherwise, we can suppose that the lowest term of the column is an isomorphism $\phi$ (by reordering the terms; note that if $X'$ contained any term of the form $\Sigma^m\Gamma_{\ell}$, we could not suppose this, because by reordering the terms, the differential of $X'$ could then not be triangular anymore).  In this case, the morphism is a section, whose retraction is given by the matrix $(0,0,\ldots,\phi^{-1})$.  Thus $X$ is quasi-isomorphic to a summand of $X'$, and is thus in $\cC^{[a,b-1]}$.  By induction hypothesis, it is quasi-isomorphic to a minimal perfect.

If $r-q+1$ is greater than one, then $X$ is obtained from $X'$ in the following recursive fashion.  Put $X_0 = X'$, and for an integer $k > 0$, let $X_k$ be the cone of the morphism
\begin{displaymath}
  \left( \begin{array}{c}
    \Sigma^{-1}f_{1,q+k-1} \\
    \Sigma^{-1}f_{2,q+k-1} \\
    \vdots  \\
    \Sigma^{-1}f_{q+k-2,q+k-1} 
  \end{array} \right)
\end{displaymath}
into $X_{k-1}$.  Then $X$ is equal to $X_{r-q+1}$.

If one of these columns contains an isomorphism, we can reorder the terms so that the isomorphism is contained in the first of these columns.  Then, by the above reasoning, this first column is a section, $X_1$ is quasi-isomorphic to a dg module which has no summands of the form $\Sigma^m\Gamma_{\ell}$, and $X$ has only $r-q$ summands of this form.  By induction, $X$ is quasi-isomorphic to a minimal perfect dg module.  This finishes the proof.

}

\demo{ (of Lemma \ref{rema::fund}.) It is easily seen that $\pr{\cD\Gamma}\Gamma$ in contained in $\cD_{\leq 0} \cap {}^{\perp}\cD_{\leq -2} \cap \perf \Gamma$.  Let $X$ be in $\cD_{\leq 0} \cap {}^{\perp}\cD_{\leq -2} \cap \perf \Gamma$.  Then $X$ is quasi-isomorphic to a minimal perfect dg module.  Thus suppose that $X$ is minimal perfect.

Let $S_i$ be the simple dg module at the vertex $i$.  Since $X$ is minimal perfect, the dimension of $\Hom{\cD\Gamma}(X, \Sigma^p S_i)$ is equal to the number of summands of $X$ isomorphic to $\Sigma^p\Gamma_i$, as a graded $\Gamma$-module.  Since $X$ is in $\cD_{\leq 0} \cap {}^{\perp}\cD_{\leq -2}$, this number is zero unless $i$ is $0$ or $1$.  This proves that $X$ is the cone of a morphism between objects of $\add\Gamma$, and thus $X$ is in $\pr{\cD\Gamma}\Gamma$.
}

We will need a particular result on the calculus of fractions in $\cC$ for certain objects.  Recall that, for any two objects $X$ and $Y$ of $\perf\Gamma$, the space $\Hom{\cC}(X,\Sigma Y)$ is the colimit of the  direct system $(\Hom{\cD\Gamma}(X', \Sigma Y))$ taken over all morphisms $f:X'\rightarrow X$ whose cone is in $\cD_{fd}\Gamma$ .

\begin{lemma}\label{lemm::fractions}
Let $X$ and $Y$ be objects of $\pr{\cD\Gamma}\Gamma$.  Then the space $\Hom{\cC}(X,\Sigma Y)$ is the colimit of the  direct system $(\Hom{\cD\Gamma}(X', \Sigma Y))$ taken over all morphisms $f:X'\rightarrow X$ whose cone is in $\cD_{fd}\Gamma \cap \cD_{\leq 0} \cap \cD_{\geq 0}$ and such that $X'$ lies in $\cD_{\leq 0}$.
\end{lemma}

\demo{There is a natural map
\begin{displaymath}
	\colim\Hom{\cD\Gamma}(X', \Sigma Y) \longrightarrow \Hom{\cC}(X, \Sigma Y),
\end{displaymath}
where the colimit is taken over all morphisms $f:X'\rightarrow X$ whose cone is in $\cD_{fd}\Gamma \cap \cD_{\leq 0} \cap \cD_{\geq 0}$ and such that $X'$ lies in $\cD_{\leq 0}$.

We first prove that it is surjective.  Let $\xymatrix{X & X'\ar[l]_{s}\ar[r]^{f} & \Sigma Y}$ be a morphism in $\cC$, with $N = cone(s)$ in $\cD_{fd}\Gamma$.

Using the canonical morphism $N \rightarrow \tau_{\geq 0}N$, we get a commuting diagram whose two lower rows and two leftmost columns are triangles:
\begin{displaymath}
	\xymatrix{ \Sigma^{-1} \tau_{< 0}N\ar@{=}[r]\ar[d] & \Sigma^{-1} \tau_{< 0}N\ar[d] & & \\
	           \Sigma^{-1} N \ar[d]\ar[r] & X' \ar[r]^s\ar[d]^a & X \ar[r]\ar@{=}[d] & N \ar[d]   \\
	           \Sigma^{-1} \tau_{\geq 0}N\ar[r]\ar[d] & X'' \ar[d]\ar[r]^t & X \ar[r] & \tau_{\geq 0}N.  \\
	           \tau_{< 0}N\ar@{=}[r] & \tau_{< 0}N & &
	}
\end{displaymath}

Thanks to the $3$-Calabi--Yau property, $\Hom{\cD\Gamma}(\Sigma^{-1}\tau_{<0}N, \Sigma Y)$ is isomorphic  to $D\Hom{\cD\Gamma}(Y, \Sigma \tau_{<0}N)$, and this is zero since $\tau_{<0}N$ is in $\cD_{\leq -2}$.  Therefore $f$ factors through $a$, and there exists a morphism $g:X''\longrightarrow \Sigma Y$ such that $ga = f$.  The fraction $\xymatrix{X & X''\ar[l]_{t}\ar[r]^{g} & \Sigma Y}$ is equal to $\xymatrix{X & X'\ar[l]_{s}\ar[r]^{f} & \Sigma Y}$, and the cone of $t$ is in $\cD_{fd}\cap \cD_{\geq 0}$.

Using the canonical morphism $\tau_{\leq 0}\tau_{\geq 0}N \rightarrow \tau_{\geq 0}N$, we get a commuting diagram whose rows are triangles:
\begin{displaymath}
	\xymatrix{ \Sigma^{-1} \tau_{\leq 0}\tau_{\geq 0}N \ar[d]\ar[r] & X''' \ar[r]^u\ar[d]^b & X \ar[r]\ar@{=}[d] & \tau_{\leq 0}\tau_{\geq 0}N \ar[d]   \\
	           \Sigma^{-1} \tau_{\geq 0}N\ar[r] & X'' \ar[r]^t & X \ar[r] & \tau_{\geq 0}N.  \\
	}
\end{displaymath}

Taking $h = bg$, we get a fraction $\xymatrix{X & X'''\ar[l]_{u}\ar[r]^{h} & \Sigma Y}$ which is equal to $\xymatrix{X & X''\ar[l]_{t}\ar[r]^{g} & \Sigma Y}$ and is such that the cone of $u$ lies in $\cD_{fd}\cap \cD_{\geq 0} \cap \cD_{\leq 0}$.

However, $X'''$ has no reason to lie in $\cD_{\leq 0}$.  Using the canonical morphism $\tau_{\leq 0}X''' \rightarrow X'''$, we get another commuting diagram whose middle rows and leftmost columns are triangles:
\begin{displaymath}
	\xymatrix{  \Sigma^{-1} \tau_{>0}X'''\ar[d]\ar@{=}[r] & \Sigma^{-1} \tau_{>0}X'''\ar[d] & & \\
	           \Sigma^{-1} M \ar[d]\ar[r] & \tau_{\leq 0}X''' \ar[r]^v\ar[d]^c & X \ar[r]\ar@{=}[d] & M \ar[d]   \\
	           \Sigma^{-1} \tau_{\leq 0}\tau_{\geq 0}N \ar[r]\ar[d] & X''' \ar[r]^u\ar[d] & X \ar[r] & \tau_{\leq 0}\tau_{\geq 0}N.  \\
	           \tau_{>0}X'''\ar@{=}[r] & \tau_{>0}X''' & & 
	}
\end{displaymath}

Since $X$ and $\tau_{\leq 0}X'''$ are in $\cD_{\leq 0}$, then so is $M$. Moreover, $\tau_{>0}X''' = \Ho^1X'''$ is in $\cD_{fd}$; indeed, the lower triangle gives an exact sequence $\Ho^0\tau_{\leq 0}\tau_{\geq 0}N \rightarrow \Ho^1X''' \rightarrow \Ho^1X$ whose leftmost term is finite-dimensional and whose rightmost term is zero. Therefore, since $\tau_{>0}X'''$ and $\tau_{\leq 0}\tau_{\geq 0}N$ are in $\cD_{\geq 0} \cap \cD_{fd}$, then so is $M$, thanks to the leftmost triangle.   

Hence, if we put $j = hc$, we have a new fraction $\xymatrix{X & \tau_{\leq 0}X'''\ar[l]_{v}\ar[r]^{j} & \Sigma Y}$ which is equal to $\xymatrix{X & X'''\ar[l]_{u}\ar[r]^{h} & \Sigma Y}$, and which is such that $\tau_{\leq 0}X'''$ is in $\cD_{\leq 0}$ and $cone(v)$ is in $\cD_{fd} \cap \cD_{\leq 0} \cap \cD_{\geq 0}$.  This proves surjectivity of the map.

We now prove that the map is injective.  Let $\xymatrix{X & X'\ar[l]_{s}\ar[r]^{f} & \Sigma Y}$ be a fraction with $X'$ in $\cD_{\leq 0}$ and $cone(s)$ in $\cD_{fd} \cap \cD_{\leq 0} \cap \cD_{\geq 0}$.  Suppose it is zero in $\Hom{\cC}(X, \Sigma Y)$, that is, $f$ factors through an object of $\cD_{fd}$.  We must prove that it factors through an object of $\cD_{fd} \cap \cD_{\leq 0} \cap \cD_{\geq 0}$.

Put $f = hg$, with $g:X'\rightarrow M$ and $h:M\rightarrow \Sigma Y$, and $M$ an object of $\cD_{fd}$.  Consider the following diagram:
\begin{displaymath}
	\xymatrix{ X'\ar@{.>}[rr]^{\varphi}\ar[rd]^{g} & & \tau_{\leq 0}\tau_{\geq 0}M\ar[d]_c \\
	           \tau_{<0}M \ar[r]^a & M\ar[r]^b\ar[dl]_{h} & \tau_{\geq 0}M\ar[d]_d\ar@{.>}[lld] \\
	           \Sigma Y & & \tau_{>0}M.
	}
\end{displaymath} 

By the $3$-Calabi--Yau property, we have an isomorphism $\Hom{\cD\Gamma}(\tau_{<0}M, \Sigma Y) \cong D\Hom{\cD\Gamma}(Y, \Sigma^2 \tau_{<0}M)$, and this is zero since $\Sigma^2\tau_{<0}M$ is in $\cD_{\leq -3}$.  Hence $h$ factors through $b$.

Moreover, $\Hom{\cD\Gamma}(X', \tau_{> 0}M)$ is zero, since $X'$ is in $\cD_{\leq 0}$ and $\tau_{>0}M$ is in $\cD_{>0}$.  Hence $bg$ factors through $c$.

This shows that $f = hg$ factors through $\tau_{\leq 0}\tau_{\geq 0}M$, which is an object of $\cD_{fd} \cap \cD_{\leq 0} \cap \cD_{\geq 0}$.  Embed $\varphi$ in a triangle
\begin{displaymath}
	\xymatrix{ X'' \ar[r]^{\varepsilon} & X'\ar[r]^{\varphi\phantom{xxx}} & \tau_{\leq 0}\tau_{\geq 0} M\ar[r] & \Sigma X''.
	}
\end{displaymath}

Then the fraction $(s\varepsilon)^{-1}(f\varepsilon)$ is equal to $s^{-1}f$.  Since $f$ factors through $\tau_{\leq 0}\tau_{\geq 0} M$, $f\varepsilon$ is zero.  

Consider finally the natural morphism $\sigma:\tau_{\leq 0}X'' \rightarrow X''$.  Its cone $\tau_{>0}X''$ is isomorphic to $\Sigma^{-1}\tau_{\leq 0}\tau_{\geq 0}M$, and is thus in $\cD_{fd}$.  Therefore the cone of $s\varepsilon\sigma$ is also in $\cD_{fd}$ by composition, and we have a fraction $(s\varepsilon\sigma)^{-1}(f\varepsilon\sigma)$ which is equal to $s^{-1}f$, and is such that $f\varepsilon\sigma = 0$, $\tau_{\leq 0}X'' \in \cD_{\leq 0}$ and $cone(s\varepsilon\sigma) \in \cD_{\leq 0} \cap \cD_{\geq 0} \cap \cD_{fd}$.  This proves injectivity of the map. 
}

Using the isomorphism of Theorem \ref{theo::3CY}, we get a bifunctorial non-degenerate bilinear form
\begin{displaymath}
	\beta_{M,L} : \Hom{\cD\Gamma}(M,L) \times \Hom{\cD\Gamma}(\Sigma^{-3}L, M) \longrightarrow k
\end{displaymath}
for $M$ in $\cD_{fd}$ and $L$ in $\perf \Gamma$.  Using this, C. Amiot constructs in \cite[Section 1.1]{Amiot08} a bifunctorial bilinear form
\begin{displaymath}
	\overline{\beta}_{X,Y} : \Hom{\cD\Gamma}(X,Y) \times \Hom{\cD\Gamma}(Y, \Sigma^2 X) \longrightarrow k
\end{displaymath}
for $X$ and $Y$ in $\cC$ in the following way.

Using the calculus of left fractions, let $s^{-1}f : X\rightarrow Y$ and $t^{-1}g : Y\rightarrow \Sigma^2 X$ be morphisms in $\cC$.  Composing them, we get a diagram
\begin{displaymath}
	\xymatrix{ X\ar[dr]^f & & Y\ar[dl]_s\ar[dr]^g & & \Sigma^2 X\ar[dl]_t \\
	                      & Y'\ar[dr]^h &   & \Sigma^2 X'\ar[dl]_{s'} & \\
	                      & & \Sigma^2 X''. & & 
	}
\end{displaymath}

Put $\Sigma^2 u = s't$.  Then one gets a commuting diagram, where rows are triangles:
\begin{displaymath}
	\xymatrix{ N\ar[r]^a & X\ar[r]^u\ar[d]^f & X''\ar[r] & \Sigma N \\
	                     & Y\ar[d]^h & & \\
	           \Sigma^2 X'\ar[r] & \Sigma^2 X''\ar[r]^b & \Sigma^2 N\ar[r] & \Sigma^3 X'.
	}
\end{displaymath}
Note that $N$ is in $\cD_{fd}$.  We put $\overline{\beta}_{X,Y}(s^{-1}f, t^{-1}g) = \beta_{N, Y'}(fa, bh)$.

\begin{proposition}\label{prop::2CY}
Let $X$ be an object of $\pr{\cC}\Gamma \cup \pr{\cC}\Sigma^{-1}\Gamma$ and $Y$ be an object of $\pr{\cC}\Gamma$.  Then the bifunctorial bilinear form
  \begin{displaymath}
	  \overline{\beta}_{X,Y} : \Hom{\cC}(X,Y) \times \Hom{\cC}(Y, \Sigma^2 X) \longrightarrow k.
  \end{displaymath}
is non-degenerate.  In particular, if one of the two spaces is finite-dimensional, then so is the other.
\end{proposition}
\demo{Let $X$ and $Y$ be objects in $\pr{\cC}\Gamma \cup \pr{\cC}\Sigma^{-1}\Gamma$ and in $\pr{\cC}\Gamma$, respectively.  In view of Proposition \ref{prop::fund}, there exist lifts $\overline{X}$ and $\overline{Y}$ of $X$ and $Y$ in $\pr{\cD\Gamma}\Gamma \cup \pr{\cD\Gamma}\Sigma^{-1}\Gamma$ and $\pr{\cD\Gamma}\Gamma$, respectively.  In particular, $\overline{X}$ and $\overline{Y}$ lie in $\cD_{\leq 1}\Gamma$.

Using the calculus of right-fractions, let $f\circ s^{-1}$ be non-zero a morphism from $X$ to $Y$ in $\cC = \perf\Gamma / \cD_{fd}\Gamma$, with $f:\overline{X}'\longrightarrow \overline{Y}$ and $s:\overline{X}'\longrightarrow \overline{X}$ morphisms in $\cD\Gamma$ such that the cone of $s$ is in $\cD_{fd}\Gamma$.  

If $\overline{X}$ lies in $\pr{\cD\Gamma}\Gamma$, then Lemma \ref{lemm::fund} allows us to suppose that $\overline{X}' = \overline{X}$ and $s = id_{\overline{X}}$.  If $\overline{X}$ lies in $\pr{\cD\Gamma}\Sigma^{-1}\Gamma$, then Lemma \ref{lemm::fractions} allows us to suppose that $\overline{X}'$ lies in $\cD_{\leq 1}$.  In both case, $\overline{X}'$ lies in $\cD_{\leq 1}$.

We now use \cite[Proposition 2.19]{KY09} : the (contravariant) functor
\begin{displaymath}
\begin{array}{cccc}
	\Phi : & \perf \Gamma & \longrightarrow & \Mod(\cD_{fd}(\Gamma)^{op}) \\	         
	       &  P           & \longmapsto     &  \Hom{\cD\Gamma}(P, ?)|_{\cD_{fd}\Gamma}
\end{array}
\end{displaymath}
is fully faithful.  Thus $\Phi(f) \neq 0$, meaning there exist $N$ in $\cD_{fd}\Gamma$ and a morphism $h:\overline{Y}\longrightarrow N$ such that its composition with $f$ is non-zero.

Recall from Theorem \ref{theo::3CY} that we have a non-degenerate bilinear form 
\begin{displaymath}
\beta : \Hom{\cD\Gamma}(\Sigma^{-3}N, X) \times \Hom{\cD\Gamma}(X,N) \longrightarrow k,
\end{displaymath} 
so there exists a morphism $j:\Sigma^{-3}N \longrightarrow X$ such that $\beta(j, h\circ f) \neq 0$.  

All the morphisms can be arranged in the following commuting diagram, where the upper and lower row are triangles in $\cD\Gamma$.
\begin{displaymath}
	\xymatrix{ \Sigma^{-3}N\ar[r]^{j} & \overline{X}'\ar[r]^{g}\ar[d]_{f} & X''\ar[r] & \Sigma^{-2}N \\
	                        & \overline{Y}\ar[r]^{h}\ar@{.>}[d]_{\ell} & N\ar@{=}[d] &              \\
	           \Sigma^{2}\overline{X}'\ar[r]^{\Sigma^2 g} & \Sigma^{2}X''\ar[r] & N\ar[r]^{\Sigma^{3}j} & \Sigma^{3}\overline{X}'.
	}
\end{displaymath}

We will show the existence of a morphism $\ell : \overline{Y} \longrightarrow \Sigma^2 X''$ making the above diagram commute.  Once this is shown, the construction of \cite{Amiot08} gives that
\begin{displaymath}
\overline{\beta}_{X,Y}(fs^{-1}, (\Sigma^2 s)(\Sigma^2 g)^{-1}\circ \ell) \ = \ \overline{\beta}_{X,Y}(f, (\Sigma^2 g)^{-1}\circ \ell) \ = \ \beta(j, h\circ f) \ \neq \ 0,
\end{displaymath}
and shows that $\overline{\beta}_{X,Y}$ is non-degenerate (here the first equality follows from the bifunctoriality of the bilinear form, and the second follows from its definition).

The existence of $\ell$ follows from the fact that $\Hom{\cD\Gamma}(\overline{Y}, \cD_{\leq -2}\Gamma) = 0$, so that $(\Sigma^3 j)\circ h = 0$.  
}

To end this section, we will prove that, in general, $\pr{\cC}\Gamma$ is not equal to the whole cluster category.

\begin{lemma}
  Let $(Q,W)$ be a quiver with potential which is not Jacobi-finite.  Then $\Sigma^2 \Gamma$ is not in $\pr{\cC}\Gamma$.
\end{lemma}
\demo{  Suppose that $\Sigma^2 \Gamma$ lies in $\pr{\cC}\Gamma$.  Then, by Proposition \ref{prop::fund}, it lifts to an object $X$ in $\pr{\cD\Gamma}\Gamma$.  We have that
\begin{displaymath}
	\Hom{\cC}(\Gamma, \Sigma^2\Gamma) = \Hom{\cD\Gamma}(\Gamma, X) = \Ho^0X.
\end{displaymath}

Now, since $X$ and $\Sigma^2\Gamma$ have the same image in $\cC$, and since $\Ho^0\Sigma^2\Gamma$ is zero (and thus finite-dimensional), $\Ho^0X$ must be finite-dimensional.

By Proposition \ref{prop::2CY}, this implies that $\Hom{\cC}(\Gamma, \Gamma)$ is also finite-dimensional, contradicting the hypothesis that $\Ho^0\Gamma = J(Q,W)$ is of infinite dimension.

Thus $\Sigma^2\Gamma$ cannot be in $\pr{\cC}\Gamma$.
}

\subsection{Mutation of $\Gamma$ in $\pr{\cD\Gamma}\Gamma$}\label{sect::nagao}

Recall the equivalence $F$ of Theorem \ref{theo::KY}.  For convenience, we shall denote it by $F^{+}$.  Denote by $F^{-}$ the quasi-inverse of the functor $\cD\Gamma \longrightarrow \cD\Gamma'$ obtained by applying Theorem \ref{theo::KY} to the mutation of $\mu_i(Q,W)$ at the vertex $i$.  Then $F^{-}(\Gamma'_i)$ is isomorphic to the cone of the morphism
\begin{displaymath}
	     \Sigma^{-1}\bigoplus_{\alpha}\Gamma_{s(\alpha)} \longrightarrow \Sigma^{-1}\Gamma_i
     \end{displaymath}
     whose components are given by right multiplication by $\alpha$.

In this subsection, we prove the following theorem, which was first formulated
as a ``hope'' by K. Nagao in his message \cite{Nag} and which is used in \cite{Nag09}.  In a more restrictive setup, an analogous result was obtained in \cite[Corollary 5.7]{IR07}.

\begin{theorem}\label{theo::mutation}
Let $\Gamma$ be the complete Ginzburg dg algebra of a quiver with potential $(Q,W)$. Let $(\varepsilon_1, \varepsilon_2, \ldots, \varepsilon_{r-1})$ be a sequence of signs.  Let $(i_1, \ldots, i_r)$ be an admissible sequence of vertices, and  let $T = \bigoplus_{j\in Q_0} T_j$ be the image of $\Gamma^{(r)}$ by the sequence of equivalences
\begin{displaymath}
	\xymatrix{\cD\Gamma^{(r)} \ar[r]^{F_{r-1}^{\varepsilon_{r-1}}} & \ldots\ar[r]^{F_1^{\varepsilon_1}\phantom{xxx}} & \cD\Gamma^{(1)} = \cD\Gamma.
	}
\end{displaymath}
Suppose that $T_j$ lies in $\pr{\cD\Gamma}\Gamma$ for all vertices $j$ of $Q$.  Then there exists a sign $\varepsilon_r$ such that all summands of the image of $\Gamma^{(r+1)}$ by $F_1^{\varepsilon_{1}}F_{2}^{\varepsilon_{2}}\cdots F_r^{\varepsilon_{r}}$ lie in $\pr{\cD\Gamma}\Gamma$.  
\end{theorem}

We start by proving a result relating morphisms in the cluster category and in the derived category, first proved in \cite[Proposition 2.12]{Amiot08} in the $\Hom{}$-finite case.

\begin{proposition}\label{prop::ext}
Let $X$ and $Y$ be objects of $\pr{\cD\Gamma}\Gamma$ such that $\Hom{\cD\Gamma}(X,\Sigma Y)$ is finite-dimensional.  Then there is an exact sequence of vector spaces
\begin{displaymath}
	0 \longrightarrow \Hom{\cD\Gamma}(X,\Sigma Y) \longrightarrow \Hom{\cC}(X,\Sigma Y) \longrightarrow D\Hom{\cD\Gamma}(Y,\Sigma X) \longrightarrow 0.
\end{displaymath}
\end{proposition}

The proof of the proposition requires some preparation.  First a lemma on limits.

\begin{lemma}\label{lemm::limits}
Let $(V_i)$ be an inverse system of finite-dimensional vector spaces with finite-dimensional limit.  Then the canonical arrow
\begin{displaymath}
	\colim (DV_i) \longrightarrow D(\Lim V_i) 
\end{displaymath}
is an isomorphism.
\end{lemma}
\demo{ This follows by duality from the isomorphisms
\begin{displaymath}
	D\colim (D V_i) \cong \Lim (DDV_i) \cong \Lim V_i.
\end{displaymath}
}

We can now prove Proposition \ref{prop::ext}.

\demo{(of Proposition \ref{prop::ext}.)  Let $X' \longrightarrow X \longrightarrow N \longrightarrow \Sigma X'$ be a triangle in $\cD\Gamma$, with $X'$ in $\cD_{\leq 0}$ and $N$ in $\cD_{fd} \cap \cD_{\leq 0}\cap\cD_{\geq 0}$.  

By the $3$-Calabi--Yau property, $\Hom{\cD\Gamma}(N, \Sigma Y) \cong D\Hom{\cD\Gamma}(Y, \Sigma^2 N)$, and this is zero since $\Sigma^2 N$ is in $\cD_{\leq -2}$.  Moreover, $\Hom{\cD\Gamma}(\Sigma^{-1}X, \Sigma Y) \cong \Hom{\cD\Gamma}(X, \Sigma^2 Y)$, and this is also zero since $\Sigma^2 Y$ is in $\cD_{\leq -2}$.

The above triangle thus gives an exact sequence
\begin{displaymath}
	0 \longrightarrow \Hom{\cD\Gamma}(X, \Sigma Y) \longrightarrow \Hom{\cD\Gamma}(X', \Sigma Y) \longrightarrow \Hom{\cD\Gamma}(\Sigma^{-1}N, \Sigma Y) \longrightarrow 0.
\end{displaymath}

We want to take the colimit of this exact sequence with respect to all morphisms $f'':X''\rightarrow X$ whose cone is in $\cD_{fd} \cap \cD_{\leq 0}\cap\cD_{\geq 0}$ and with $X''$ in $\cD_{\leq 0}$.  The colimit will still be a short exact sequence, since, as we shall prove, all the spaces involved and their colimits are finite-dimensional.

The leftmost term is constant; its colimit is itself.

Consider the rightmost term.  Since $Y$ is in $\pr{\cD\Gamma}\Gamma$, there is a triangle
\begin{displaymath}
	P_1 \longrightarrow P_0 \longrightarrow Y \longrightarrow \Sigma P_1
\end{displaymath}
with $P_0$ and $P_1$ in $\add \Gamma$.  Noticing that $\Hom{\cD\Gamma}(P_i, N) = \Hom{\cD\Gamma}(\Ho^0P_i, \Ho^0N)$ for $i \in \{1,2\}$, this yields get an exact sequence
\begin{displaymath}
	0 \rightarrow \Hom{\cD}(Y,N) \rightarrow \Hom{\cD}(\Ho^0P_0,\Ho^0N) \rightarrow \Hom{\cD}(\Ho^0P_1,\Ho^0N) \rightarrow \Hom{\cD}(Y, \Sigma N) \rightarrow 0.
\end{displaymath}

Since the two middle spaces are finite-dimensional, so are the other two, and the limit of this sequence is still exact.  

Now $H^0X$ is an $(\End{\cD\Gamma}\Gamma)$-module.  Since $\End{\cD\Gamma}\Gamma$ is the jacobian algebra of a quiver with potential, $H^0X$ is the limit of all its finite-dimensional quotients.  The system given by the $H^0N$ is a system of all the finite-dimensional quotients of $H^0X$; its limit is thus $H^0X$. 

 Hence the limit of $\Hom{\cD\Gamma}(\Ho^0P_i, \Ho^0 N)$ is $\Hom{\cD\Gamma}(\Ho^0P_i, \Ho^0X)$, which is isomorphic to $\Hom{\cD\Gamma}(P_i, \Ho^0X)$. We thus have an exact sequence 
\begin{displaymath}
	\Hom{\cD\Gamma}(P_0,\Ho^0X) \rightarrow \Hom{\cD\Gamma}(P_1,\Ho^0X) \rightarrow \Lim \Hom{\cD\Gamma}(Y, \Sigma N) \rightarrow 0.
\end{displaymath} 

This implies the isomorphisms
 
\begin{displaymath}
	\Lim \Hom{\cD\Gamma}(Y, \Sigma N) \cong \Hom{\cD\Gamma}(Y, \Sigma\Ho^0 X) \cong \Hom{\cD\Gamma}(Y, \Sigma X).
\end{displaymath} 

Using Lemma \ref{lemm::limits}, we thus get that $\colim D\Hom{\cD\Gamma}(Y, \Sigma N) = D\Hom{\cD\Gamma}(Y, \Sigma X)$, and the $3$-Calabi--Yau property of Theorem \ref{theo::3CY} implies that $D\Hom{\cD\Gamma}(Y, \Sigma N)$ is isomorphic to $\Hom{\cD\Gamma}(\Sigma^{-1}N, \Sigma Y)$.  Therefore the colimit of the $\Hom{\cD\Gamma}(\Sigma^{-1}N, \Sigma Y)$ is $D\Hom{\cD\Gamma}(Y, \Sigma X)$ as desired.

It remains to be shown that the colimit of the terms of the form $\Hom{\cD\Gamma}(X', \Sigma Y)$ is $\Hom{\cD\Gamma}(X, \Sigma Y)$.  This is exactly Lemma \ref{lemm::fractions}.  This finishes the proof of the Proposition.
}

This enables us to formulate a result on the lifting of triangles from the cluster category to the derived category.

\begin{proposition}\label{prop::lifting}
Let $X$ and $Y$ be objects of $\pr{\cC}\Gamma$, with $\dim \Hom{\cC}(X, \Sigma Y) = 1$ (and so $\dim \Hom{\cC}(Y, \Sigma X) = 1$ by Proposition \ref{prop::2CY}).  Let
 \begin{displaymath}
   X \longrightarrow E \longrightarrow Y \longrightarrow \Sigma X \ \textrm{and} \	Y \longrightarrow E' \longrightarrow X \longrightarrow \Sigma Y
 \end{displaymath}
be non-split triangles (they are unique up to isomorphism).  Then one of the two triangles lifts to a triangle $A \longrightarrow B \longrightarrow C \longrightarrow \Sigma A$ in $\perf \Gamma$, with $A$, $B$ and $C$ in $\pr{\cD\Gamma}\Gamma$.  
\end{proposition}
\demo{ According to Proposition \ref{prop::fund}, we can lift $X$ and $Y$ to objects $\overline{X}$ and $\overline{Y}$ of $\pr{\cD}\Gamma$.  Using the short exact sequence of Proposition \ref{prop::ext}, we have that one of $\Hom{\cD\Gamma}(\overline{X}, \Sigma\overline{Y})$ and $\Hom{\cD\Gamma}(\overline{Y}, \Sigma\overline{X})$ is one-dimensional.

Suppose that $\Hom{\cD\Gamma}(\overline{Y}, \Sigma\overline{X})$ is one-dimensional.  Let
\begin{displaymath}
	\overline{X} \longrightarrow \overline{E} \longrightarrow \overline{Y} \longrightarrow \Sigma\overline{X}
\end{displaymath}
be a non-split triangle.  Since $\pr{\cD\Gamma}\Gamma = \cD_{\leq 0} \cap {}^{\perp}\cD_{\leq -2} \cap \perf\Gamma$ is closed under extensions, $\overline{E}$ lies in $\pr{\cD\Gamma}\Gamma$.  Thus the equivalence of Proposition \ref{prop::fund} implies that the triangle descends to a non-split triangle in $\cC$.  Up to isomorphism, this non-split triangle is $X \longrightarrow E \longrightarrow Y \longrightarrow \Sigma X$.

The proof is similar if $\Hom{\cD\Gamma}(\overline{X}, \Sigma\overline{Y})$ is one-dimensional; in this case, the triangle $Y \longrightarrow E' \longrightarrow X \longrightarrow \Sigma Y$ is the one which can be lifted.
}

We can now prove the main theorem of this subsection.

\demo{ (of Theorem \ref{theo::mutation}.)

Put $i = i_r$.  For any vertex $j \neq i$, the image of $\Gamma^{(r+1)}$ by $F_1^{\varepsilon_{1}}F_{2}^{\varepsilon_{2}}\cdots F_{r-1}^{\varepsilon_{r-1}}F_r^{\varepsilon}$ is isomorphic to $T_j$ for any sign $\varepsilon$, and is in $\pr{\cD\Gamma}\Gamma$ by hypothesis.  Now, the images of $\Gamma^{(r+1)}_i$ by $F_1^{\varepsilon_{1}}F_{2}^{\varepsilon_{2}}\cdots F_{r-1}^{\varepsilon_{r-1}}F_r^{+}$ and by $F_1^{\varepsilon_{1}}F_{2}^{\varepsilon_{2}}\cdots F_{r-1}^{\varepsilon_{r-1}}F_r^{-}$ become isomorphic in the cluster category $\cC_{Q,W}$, and they lie in $\pr{\cC}\Gamma$.  Denote these images by $T_i^*$.  We have that $\dim \Hom{\cC}(T_i, \Sigma T_i^*) = 1$.  Thus we can apply Proposition \ref{prop::lifting} and get that $T_i^*$ is lifted in $\pr{\cD\Gamma}\Gamma$ either to $F_1^{\varepsilon_{1}}F_{2}^{\varepsilon_{2}}\cdots F_{r-1}^{\varepsilon_{r-1}}F_r^{+}(\Gamma^{r+1}_i)$ or to $F_1^{\varepsilon_{1}}F_{2}^{\varepsilon_{2}}\cdots F_{r-1}^{\varepsilon_{r-1}} F_r^{-}(\Gamma^{r+1}_i)$.
}

\section{Cluster character}\label{sect::clusterchar}

Let $\cC$ be a (not necessarily Hom--finite) triangulated category with suspension functor $\Sigma$.  Let $T = \bigoplus_{i=1}^{n}T_i$ be a basic rigid object in $\cC$ (with each $T_i$ indecomposable), that is, an object $T$ such that $\Hom{\cC}(T, \Sigma T) = 0$ and $i \neq j$ implies that $T_i$ and $T_j$ are not isomorphic.  We will assume the following :

\begin{enumerate}
	\item $\pr{\cC}T$ is a Krull--Schmidt category;
	\item $B = \End{\cC}(T)$ is the (completed) Jacobian algebra of a quiver with potential $(Q,W)$;
	\item the simple $B$--module at each vertex can be lifted to an object in $\pr{\cC}(T) \cap \pr{\cC}(\Sigma T)$ through the functor $\Hom{\cC}(T,-)$;
	\item for all objects $X$ of $\pr{\cC}(\Sigma T) \cup \pr{\cC}(T)$ and $Y$ of $\pr{\cC}(\Sigma T)$, there exists a non-degenerate bilinear form
	  \begin{displaymath}
	    \Hom{\cC}(X,Y) \times \Hom{\cC}(Y, \Sigma^2 X) \longrightarrow k
    \end{displaymath}
which is functorial in both variables.
\end{enumerate}

\begin{lemma}
The above hypotheses hold for the cluster category $\cC_{Q,W}$ of a quiver with potential $(Q,W)$, where $T$ is taken to be $\Sigma^{-1}\Gamma$.
\end{lemma}
\demo{ Condition (1) is proved in Corollary \ref{coro::prkrull}, since $\pr{\cC}\Gamma$ is equivalent to $\pr{\cC}\Sigma^{-1}\Gamma$.
Condition (2) follows from Proposition \ref{prop::fund}, since $\End{\cC}(\Sigma^{-1}\Gamma)$ is isomorphic to $\End{\cC}(\Gamma)$, which is in turn isomorphic to $\End{\cD\Gamma}(\Gamma) = \Ho^0\Gamma$, and this is the completed Jacobian algebra of $(Q,W)$.  Condition (3) follows from the fact that $\Hom{\cC}(\Gamma_i, \Sigma\Gamma_i^*) = \Hom{\cD\Gamma}(\Gamma_i, \Sigma\Gamma_i^*)$ is one-dimensional (see \cite[Section 4]{KY09}).  Finally, condition (4) is exactly Proposition \ref{prop::2CY}.
}

As in \cite{DK08} and \cite{Palu08}, define the \emph{index with respect to $T$} of an object $X$ of $\pr{\cC}T$ as the element of $\K_0 (\proj B)$ given by
\begin{displaymath}
	\ind{T}X = [FT_0^X] - [FT_1^X],
\end{displaymath}
where $T_1^X \longrightarrow T_0^X\longrightarrow X\longrightarrow\Sigma T_1^X$ is an $(\add T)$-presentation of $X$.  One can show as in \cite{Palu08} that the index is well-defined, that is, does not depend on the choice of a presentation.
\subsection{Modules}

Consider the functors $F = \Hom{\cC}(T,-) : \cC \longrightarrow \Mod B$ and $G = \Hom{\cC}(-,\Sigma^2 T) : \cC \longrightarrow \Mod B^{op}$, where $\Mod B$ is the category of right $B$--modules.  

For an object $U$ of $\cC$, let $(U)$ be the ideal of morphisms in $\cC$ factoring through an object of $\add U$.

This subsection is devoted to proving some useful properties of the functors $F$ and $G$.

\begin{lemma}\label{lemm::modules}
Let $X$ and $Y$ be objects in $\cC$.

\begin{enumerate}
  \item If $X$ lies in $\pr{\cC}T$, then $F$ induces an isomorphism 
    \begin{displaymath}
      \Hom{\cC}(X,Y) / (\Sigma T) \longrightarrow \Hom{B}(FX,FY).
    \end{displaymath}
    
    If $Y$ lies in $\pr{\cC}\Sigma T$, then $G$ induces an isomorphism
    \begin{displaymath}
      \Hom{\cC}(X,Y) / (\Sigma T) \longrightarrow \Hom{B^{op}}(GY,GX).
    \end{displaymath}
  
  \item $F$ induces an equivalence of categories
    \begin{displaymath}
      \pr{\cC}T / (\Sigma T) \longrightarrow \MOD B,
    \end{displaymath}
    where $\MOD B$ denotes the category of finitely presented $B$--modules.
    
  \item Any finite-dimensional $B$--module can be lifted through $F$ to an object in $\pr{\cC}T \cap \pr{\cC} \Sigma T$.  Any short exact sequence of finite-dimensional $B$-modules can be lifted through $F$ to a triangle of $\cC$, whose three terms are in $\pr{\cC}T \cap \pr{\cC}\Sigma T$.
\end{enumerate}
\end{lemma}
\demo{
(1) We only prove the first isomorphism; the proof of the second one is dual.  First, suppose that $X = T_i$ is an indecomposable summand of $T$.  Let $f:FT_i \longrightarrow FY$ be a morphism of $B$--modules.  Note that any element $g$ of $FT_i = \Hom{\cC}(T,T_i)$ is of the form $pg'$, where $p:T\longrightarrow T_i$ is the canonical projection and $g'$ is an endomorphism of $T$.  Hence $f(g) = f(p)g'$.  Moreover, consider the idempotent $e_i$ in $\End{\cC}T$ associated with $T_i$.  We have that $f(p) = f(pe_i) = f(p)e_i$.  Hence $f(p)$ can be viewed as a morphism from $T_i$ to $Y$, and $f = F(f(p))$.  This shows that there is a bijection $\Hom{\cC}(T_i,Y) \longrightarrow \Hom{B}(FT_i,FY)$.
 
One easily sees that this bijection will also hold if $X$ is a direct sum of direct summands of $T$.

Now, let $X$ be in $\pr{\cC}T$, and let $\xymatrix{T_1^X \ar[r]^{\alpha} & T_0^X \ar[r]^{\beta} & X \ar[r]^{\gamma} & \Sigma T_1^X}$ be a triangle in $\cC$, with $T_0^X, T_1^X \in \add T$.  Let $f:FX \longrightarrow FY$ be a morphism of $B$--modules.  We have that $fF\beta$ belongs to $\Hom{B}(FT_0^X, FY)$, and by the above lifts to a morphism $\omega:T_0^X \longrightarrow Y$.

Moreover, $F(\omega\alpha) = F\omega F\alpha = fF\beta F\alpha = 0$, and by injectivity, $\omega\alpha = 0$.  Hence there exists $\phi:X \longrightarrow Y$ such that $\phi\beta = \omega$, so $F\phi F\beta = f F\beta$. Since $F\beta$ is surjective, this gives $F\phi = f$.  Therefore the map $\Hom{\cC}(X,Y) \longrightarrow \Hom{B}(FX,FY)$ is surjective.

Suppose now that $u:X\longrightarrow Y$ is such that $Fu = 0$.  Then $F(u\beta) = FuF\beta = 0$, and by the injectivity proved above, $u\beta = 0$, and $u$ factors through $\Sigma T_1^X$.  This finishes the proof.

 (2) It follows from part (1) that the functor is fully faithful.  Let now $M \in \MOD B$, and let $P_1 \longrightarrow P_0 \longrightarrow M \longrightarrow 0$ be a projective presentation.  By part (1), $P_1 \longrightarrow P_0$ lifts to a morphism $T_1 \longrightarrow T_0$ in $\cC$, with $T_0, T_1 \in \add T$.  We can embed this morphism in a triangle $T_1 \longrightarrow T_0 \longrightarrow X \longrightarrow \Sigma T_1$, and we see that $FX$ is isomorphic to $M$.  This proves the equivalence.

 (3) By our hypothesis, the statement is true for the simple modules at each vertex.  Let $M$ be a finite-dimensional $B$--module.  According to a remark following Definition 10.1 of \cite{DWZ08}, $M$ is nilpotent.  Therefore it can be obtained from the simple modules by repeated extensions.  All we have to do is show that the property is preserved by extensions in $\Mod B$.

Let $0 \longrightarrow L \longrightarrow M \longrightarrow N \longrightarrow 0$ be a short exact sequence, with $L$ and $N$ in $\MOD B$ admitting lifts $\overline{L}$ and $\overline{N}$ in $\pr{\cC}T \cap \pr{\cC}\Sigma T$, respectively.  Using projective presentations of $L$ and $N$, we consctruct one for $M$ and obtain a diagram as below, where the upper two rows are split.

\begin{displaymath}
	\xymatrix{ 0 \ar[r] & P_1^L \ar[r]\ar[d] & P_1^L \oplus P_1^N \ar[r]\ar[d] & P_1^N \ar[r]\ar[d] & 0 \\
	           0 \ar[r] & P_0^L \ar[r]\ar[d] & P_0^L \oplus P_0^N \ar[r]\ar[d] & P_0^N \ar[r]\ar[d] & 0 \\
	           0 \ar[r] & L     \ar[r]\ar[d] & M                  \ar[r]\ar[d] & N     \ar[r]\ar[d] & 0 \\
	                    & 0 & 0 & 0 &
	}
\end{displaymath}

Thanks to part (2), the upper left square can be lifted into a commutative diagram

\begin{displaymath}
	\xymatrix{ T_1^L \ar[r]\ar[d] & T_1^L \oplus T_1^N \ar[d] \\
	           T_0^L \ar[r]       & T_0^L \oplus T_0^N
	}
\end{displaymath}

which in turn embeds in a nine-diagram as follows.

\begin{displaymath}
	\xymatrix{  T_1^L \ar[r]\ar[d] & T_1^L \oplus T_1^N \ar[r]\ar[d] & T_1^N \ar[r]\ar[d] & \Sigma T_1^L \\
	            T_0^L \ar[r]\ar[d] & T_0^L \oplus T_0^N \ar[r]\ar[d] & T_0^N \ar[r]\ar[d] & \Sigma T_0^L \\
	            \overline{L} \ar[r]\ar[d] & \overline{M}\ar[r]\ar[d] & \overline{N}\ar[r]\ar[d] & \Sigma\overline{L}  \\
	            \Sigma T_1^L \ar[r] & \Sigma T_1^L \oplus \Sigma T_1^N \ar[r] & \Sigma T_1^N  & 	            
	}
\end{displaymath}

Hence $\overline{M}$ is a lift of $M$ in $\pr{\cC}T$. Now, since $\overline{N}$ lies in $\pr{\cC}\Sigma T$, it follows from part (1) that the morphism $\Sigma^{-1} \overline{N} \longrightarrow \overline{L}$ is in $(\Sigma T)$, and thus from Lemma \ref{lemm::index} below that $\overline{M}$ is also in $\pr{\cC} \Sigma T$.  This finishes the proof.
}

\begin{lemma}\label{lemm::mono}
Let $X$, $Y$ and $Z$ be objects in $\cC$.  Suppose that $Y$ and $Z$ lie in $\pr{\cC}\Sigma T$ and that $FY$ is finite-dimensional. Let
\begin{displaymath}
	\xymatrix{ X \ar[r]^{f} & Y \ar[r]^{g} & Z \ar[r] & \Sigma X
	}
\end{displaymath} 
be a triangle.
If $Ff=0$, then $f \in (\Sigma T)$.
\end{lemma}
\demo{ The equality $Ff = 0$ means that $Fg$ is injective.  Using the non-degenerate bilinear form, we get a commuting diagram
\begin{displaymath}
  \xymatrix{ \Hom{\cC}(T,Y) \ar@{^{(}->}[r]^{Fg}\ar@{^{(}->}[d] & \Hom{\cC}(T,Z)\ar@{^{(}->}[d] \\
             D\Hom{\cC}(Y, \Sigma^2 T)\ar[r]^{DGg} & D\Hom{\cC}(Z, \Sigma^2 T),
  }	
\end{displaymath} 
where the top horizontal morphism and the two vertical ones are injective.  Since $FY$ is finite-dimensional, the left morphism is an isomorphism.  Thus $DGg$ is injective, and $Gg$ is surjective.  But this means that $Gf = 0$, and by part (1) of Lemma \ref{lemm::modules}, $f \in (\Sigma T)$. 
}

\subsection{Presentations and index}

Let us now study some closure properties of $\pr{\cC}T$, and deduce some relations between triangles and indices.

\begin{lemma}\label{lemm::index}
Let $\xymatrix{X \ar[r] & Y \ar[r] & Z\ar[r]^{\varepsilon} & \Sigma X}$ be a triangle in $\cC$ such that $\varepsilon$ is in $(\Sigma T)$.  Then
\begin{enumerate}
  \item If two of $X$, $Y$ and $Z$ lie in $\pr{\cC}T$, then so does the third one. 
  
  \item If $X,Y,Z \in \pr{\cC}T$, then we have an equality $\ind{T} X + \ind{T} Z = \ind{T} Y$.
\end{enumerate} 
\end{lemma}
\demo{Let us first suppose that $X$ and $Z$ lie in $\pr{\cC}T$.  Let $T_1^X\longrightarrow T_0^X\longrightarrow X \longrightarrow \Sigma T_1^X$ and $T_1^Z\longrightarrow T_0^Z\longrightarrow Z \longrightarrow \Sigma T_1^Z$ be two triangles, with $T_0^X, T_1^X, T_0^Z$ and $T_1^Z$ in $\add T$.

Since $\Hom{\cC}(T, \Sigma T) = 0$, the composition $T_0^Z\longrightarrow Z \longrightarrow \Sigma X$ vanishes, so $T_0^Z\longrightarrow Z$ factors through $Y$.  This gives a commutative square
\begin{displaymath}
	\xymatrix{T_0^X \oplus T_0^Z \ar[r]\ar[d] & T_0^Z\ar[d] \\
	                    Y        \ar[r]       & Z }
\end{displaymath} 
which can be completed into a nine-diagram
\begin{displaymath}
	\xymatrix{T_1^X \ar[r]\ar[d] & T_1^X \oplus T_1^Z \ar[r]\ar[d] & T_1^Z\ar[d]\ar[r] & \Sigma T_1^X \\
	          T_0^X \ar[r]\ar[d] & T_0^X \oplus T_0^Z \ar[r]\ar[d] & T_0^Z\ar[d]\ar[r] & \Sigma T_0^X \\
	               X\ar[r]\ar[d] &    Y        \ar[r]\ar[d]        & Z \ar[r]\ar[d]    & \Sigma X \\
	           \Sigma T_1^X      & \Sigma T_1^X \oplus \Sigma T_1^Z& \Sigma T_1^Z
}
  \end{displaymath}
showing that $Y$ is in $\pr{\cC}T$ and that assertion 2 is true.

Now suppose that $X$ and $Y$ lie in $\pr{\cC}T$.  Since the composition $Z\longrightarrow \Sigma X \longrightarrow \Sigma^2 T_1^X$ is zero, the morphism $Z\longrightarrow \Sigma X$ factors through $\Sigma T_0^X$.  This yields an octahedron

$\phantom{an octahedron}$\begin{xy} 0;<1pt,0pt>:<0pt,-1pt>:: 
(105,0) *+{\Sigma Y} ="0",
(74,90) *+{Z} ="1",
(207,90) *+{\Sigma X} ="2",
(0,69) *+{\Sigma W} ="3",
(103,147) *+{\Sigma T_0^X} ="4",
(133,69) *+{\Sigma^2 T_1^X} ="5",
"0", {\ar|+"1"},
"2", {\ar"0"},
"3", {\ar"0"},
"0", {\ar@{.>}"5"},
"1", {\ar"2"},
"3", {\ar|+"1"},
"1", {\ar"4"},
"4", {\ar"2"},
"2", {\ar@{.>}"5"},
"4", {\ar"3"},
"5", {\ar@{.>}|+"3"},
"5", {\ar@{.>}|+"4"},
\end{xy} 

which produces a triangle $T_1^X\longrightarrow W\longrightarrow Y \longrightarrow \Sigma T_1^X$.  Composing with $Y \longrightarrow \Sigma T_1^Y$, we get a second octahedron

$\phantom{an octahedron}$
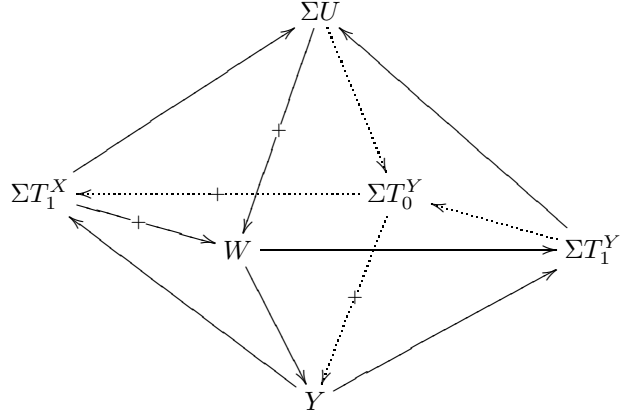
\begin{xy} 0;<1pt,0pt>:<0pt,-1pt>:: 
(105,0) *+{\Sigma U} ="0",
(74,90) *+{W} ="1",
(207,90) *+{\Sigma T_1^Y} ="2",
(0,69) *+{\Sigma T_1^X} ="3",
(103,147) *+{Y} ="4",
(133,69) *+{\Sigma T_0^Y} ="5",
"0", {\ar|+"1"},
"2", {\ar"0"},
"3", {\ar"0"},
"0", {\ar@{.>}"5"},
"1", {\ar"2"},
"3", {\ar|+"1"},
"1", {\ar"4"},
"4", {\ar"2"},
"2", {\ar@{.>}"5"},
"4", {\ar"3"},
"5", {\ar@{.>}|+"3"},
"5", {\ar@{.>}|+"4"},
\end{xy} 

which gives triangles $T_1^X\longrightarrow U\longrightarrow T_0^Y \longrightarrow \Sigma T_1^X$ and $T_1^Y\longrightarrow U\longrightarrow W \longrightarrow \Sigma T_1^Y$.  Note that, since $\Hom{\cC}(T, \Sigma T) = 0$, the first triangle is split, so $U$ is isomorphic to $T_1^X \oplus T_0^Y$.

From the first octahedron, one gets a triangle $T_0^X\longrightarrow W\longrightarrow Z \longrightarrow \Sigma T_0^X$. Construct one last octahedron with the composition $U \longrightarrow W \longrightarrow Z$.

$\phantom{an octahedron}$
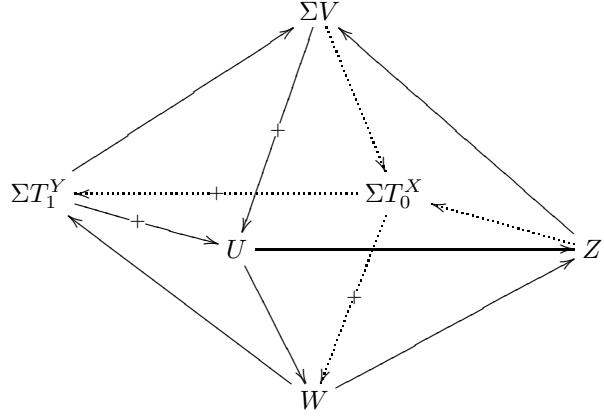
\begin{xy} 0;<1pt,0pt>:<0pt,-1pt>:: 
(105,0) *+{\Sigma V} ="0",
(74,90) *+{U} ="1",
(207,90) *+{Z} ="2",
(0,69) *+{\Sigma T_1^Y} ="3",
(103,147) *+{W} ="4",
(133,69) *+{\Sigma T_0^X} ="5",
"0", {\ar|+"1"},
"2", {\ar"0"},
"3", {\ar"0"},
"0", {\ar@{.>}"5"},
"1", {\ar"2"},
"3", {\ar|+"1"},
"1", {\ar"4"},
"4", {\ar"2"},
"2", {\ar@{.>}"5"},
"4", {\ar"3"},
"5", {\ar@{.>}|+"3"},
"5", {\ar@{.>}|+"4"},
\end{xy} 

As was the case for $U$, $V$ is in a split triangle, and is thus isomorphic to $T_1^Y \oplus T_0^X$.  Hence there is a triangle $V \longrightarrow U \longrightarrow Z \longrightarrow \Sigma V$, with $U$ and $V$ in add $T$.  This proves that $Z$ lies in $\textrm{pr}_{\mathcal{C}}T$.

Finally, suppose that $Y$ and $Z$ are in $\textrm{pr}_{\mathcal{C}}T$.  Notice that since $\Sigma^{-1} \varepsilon$ factors through $\add T$, the composition $\Sigma^{-1}T_0^Z \longrightarrow \Sigma^{-1}Z \longrightarrow X$ vanishes. Applying a reasonning dual to that of the preceding case, one proves that $X$ lies in $\pr{\cC}T$.
}
 
The next lemma is an adapted version of Proposition 6 of \cite{Palu08}.

\begin{lemma}\label{lemm::index2}
Let $\xymatrix{ X\ar[r]^{\alpha} & Y\ar[r]^{\beta} & Z\ar[r]^{\gamma} & \Sigma X }$ be a triangle in $\cC$, with $X,Z \in \pr{\cC}T$ such that $\Coker F\beta$  is finite-dimensional.  Let $C\in \pr{\cC}T \cap \pr{\cC}\Sigma T$ be such that $FC = \Coker F\beta$.  Then $Y\in \pr{\cC}T$, and $\ind{T} X + \ind{T} Z = \ind{T} Y + \ind{T} C + \ind{T} \Sigma^{-1} C$.  
\end{lemma}
\demo{ Note that since $\Coker F\beta$ is finite-dimensional, it can be lifted to $C\in \pr{\cC}T \cap \pr{\cC}\Sigma T$ thanks to Lemma \ref{lemm::modules}.

The case where $\gamma$ factors through $\add \Sigma T$ was treated in Lemma \ref{lemm::index}.  In that case, $\Coker F\beta = 0$, and $C \in \add \Sigma T$, so that $\ind{T}C = - \ind{T}\Sigma^{-1}C$.

Suppose now that $\gamma$ is not in $(\Sigma T)$.  In $\MOD B$, there is a commutative triangle

$\phantom{xxxxxxxxxxxxxx}\xymatrix{ FZ \ar[rr]\ar[dr] & & F\Sigma X \\
                       & \Coker F\beta \ar[ur] &
}$

which, thanks to Lemma \ref{lemm::modules}, we can lift to a commutative triangle

$\phantom{xxxxxxxxxxxxxxx}\xymatrix{ Z \ar[rr]^{a}\ar[dr]^{b} & & \Sigma X \\
                       & C \oplus \Sigma\overline{T} \ar[ur]^{c} &
}$

in $\cC$, where $\overline{T}$ lies in $\add T$.  Form an octahedron

$\phantom{an octahedron}$\begin{xy} 0;<1pt,0pt>:<0pt,-1pt>:: 
(105,0) *+{\Sigma Y} ="0",
(74,90) *+{Z} ="1",
(207,90) *+{\Sigma X.} ="2",
(0,69) *+{U} ="3",
(103,147) *+{C \oplus \Sigma\overline{T}} ="4",
(133,69) *+{U'} ="5",
"0", {\ar|+"1"},
"2", {\ar"0"},
"3", {\ar"0"},
"0", {\ar@{.>}"5"},
"1", {\ar"2"},
"3", {\ar|+"1"},
"1", {\ar"4"},
"4", {\ar"2"},
"2", {\ar@{.>}"5"},
"4", {\ar"3"},
"5", {\ar@{.>}|+"3"},
"5", {\ar@{.>}|+"4"},
\end{xy} 

Since $Fb$ is an epimorphism, the morphism $C \oplus \Sigma \overline{T} \longrightarrow U$ must lie in  $(\Sigma T)$, by Lemma \ref{lemm::modules}, part (1).  Since $Fc$ is a monomorphism, the same must hold for $\Sigma^{-1}U' \longrightarrow C \oplus \Sigma \overline{T}$, by Lemma \ref{lemm::mono}.  By composition, the morphism $\Sigma^{-1}U' \longrightarrow U$ is also in $(\Sigma T)$.

We thus have three triangles

$\phantom{yeahyeahye}$\xymatrix{
	\Sigma^{-1} U \ar[r]& Z \ar[r]& C \oplus \Sigma\overline{T} \ar[r]& U	\\
	\Sigma^{-1} C \oplus \overline{T} \ar[r]& X \ar[r]& \Sigma^{-1}U' \ar[r]& C\oplus \Sigma \overline{T} \\
	\Sigma^{-1} U \ar[r]& Y \ar[r]& \Sigma^{-1}U' \ar[r]& U 
}

whose third morphism factors through $\Sigma T$.  Applying Lemma \ref{lemm::index}, we get that $\Sigma^{-1} U$, $\Sigma^{-1} U'$ and $Y$ are in $\pr{\cC}T$, and that
\begin{displaymath}
	\ind{T} \Sigma^{-1}U + \ind{T} C + \ind{T} \Sigma \overline{T} = \ind{T} Z,
\end{displaymath}
\begin{displaymath}
	\ind{T} \Sigma^{-1}C + \ind{T} \overline{T} + \ind{T} \Sigma^{-1}U' = \ind{T} X,\textrm{ and}
\end{displaymath}
\begin{displaymath}
	\ind{T} Y = \ind{T} \Sigma^{-1}U + \ind{T} \Sigma^{-1}U'.
\end{displaymath}

Summing up, and noticing that $\ind{T} \overline{T} = - \ind{T} \Sigma\overline{T}$, we get the desired equality.
} 

\begin{lemma}\label{lemm::ci}
Let $X$ be an object in $\pr{\cC}T \cap \pr{\cC}\Sigma T$ such that $FX$ is finite-dimensional.  Then the sum $\ind{T}X + \ind{T}\Sigma^{-1}X$ only depends on the dimension vector of $FX$.
\end{lemma}
\demo{ First, notice that $FX = 0$ if, and only if, $X$ is in $\add \Sigma T$.

Second, suppose $X$ is indecomposable.  If $Y$ is another such object such that $FX$ and $FY$ are isomorphic and non-zero, then $X$ and $Y$ are isomorphic in $\cC$. Indeed, in view of Lemma \ref{lemm::modules}, part (2), there exist morphisms $f:X\longrightarrow Y$ and $g:Y\longrightarrow X$ such that $f\circ g = id_X + t$ and $g\circ f = id_Y + t'$, with $t, t' \in (\Sigma T)$.  But since the endomorphism rings of $X$ and $Y$ are local and not contained in $(\Sigma T)$, this implies that $f\circ g$ and $g\circ f$ are isomorphisms.

Third, let us show that the sum depends only on the isomorphism class of $FX$.  Let $Y$ be another such object such that $FX$ and $FY$ are isomorphic.  Write $X = \overline{X} \oplus \Sigma T^X$ and $Y = \overline{Y} \oplus \Sigma T^Y$, where $T^X, T^Y \in \add T$ and $\overline{X}, \overline{Y}$ have no direct summand in $\add \Sigma T$.  Then $F\overline{X} = F\overline{Y}$, and by the above $\overline{X}$ and $\overline{Y}$ are isomorphic.  We have
\begin{eqnarray*}
	\ind{T}X + \ind{T}\Sigma^{-1}X & = & \ind{T}\overline{X} + \ind{T}\Sigma^{-1}\overline{X} + \ind{T}\Sigma T^X + \ind{T}T^X\\
	 & = & \ind{T}\overline{X} + \ind{T}\Sigma^{-1}\overline{X} \\
	 & = & \ind{T}\overline{Y} + \ind{T}\Sigma^{-1}\overline{Y} + \ind{T}\Sigma T^Y + \ind{T}T^Y \\
	 & = & \ind{T}Y + \ind{T}\Sigma^{-1}Y.
\end{eqnarray*}

Finally, we prove that the sum only depends on the dimension vector of $FX$. Let $0\longrightarrow L\longrightarrow M\longrightarrow N\longrightarrow 0 $ be an exact sequence in $\MOD B$.  As in the proof of part (3) of Lemma \ref{lemm::modules}, lift it to a triangle $\overline{L}\longrightarrow \overline{M}\longrightarrow \overline{N}\longrightarrow \Sigma\overline{L}$, where the last morphism in in $(\Sigma T) \cap (\Sigma^2 T)$. Using Lemma \ref{lemm::index}, we get the equality
\begin{eqnarray*}
	\ind{T} \overline{M} + \ind{T} \Sigma^{-1}\overline{M} & = & \ind{T} \overline{L} + \ind{T} \Sigma^{-1}\overline{L} + \ind{T} \overline{N} + \ind{T} \Sigma^{-1}\overline{N}.
\end{eqnarray*}  
This gives the independance on the dimension vector.    
}

\begin{notation}\label{nota::iota}
For a dimension vector $e$, denote by $\iota(e)$ the quantity $\ind{T}X + \ind{T}\Sigma^{-1}X$, where $\dimv FX = e$ (by the above Lemma, this does not depend on the choice of such an $X$).
\end{notation}

\begin{lemma}\label{lemm::bil}
If $X \in \pr{C}(T)$ and $Y \in \pr{C}(\Sigma T)$, then the bilinear form induces a non-degenerate bilinear form
    \begin{displaymath}
	    (\Sigma T)(X,Y) \times \Hom{\cC}(Y, \Sigma^2 X)/(\Sigma^2 T) \longrightarrow k.
    \end{displaymath}
\end{lemma}
\demo{ Let $\xymatrix{ T_1^X\ar[r] & T_0^X\ar[r] & X\ar[r]^{\eta} & \Sigma T_1^X }$ be a triangle, with $T_0^X, T_1^X$ in $\add T$.  Consider the following diagram.
\begin{displaymath}
	\xymatrix{ \Hom{C}(X,Y) & \times & \Hom{C}(Y, \Sigma^2 X)\ar[d]^{\Sigma^2 \eta_*}  & \ar[r]& k \\
	           \Hom{C}(\Sigma T_1^X,Y)\ar[u]_{\eta^*} & \times & \Hom{C}(Y, \Sigma^3 T_1^X)  & \ar[r]& k
	}
\end{displaymath}

The bifunctoriality of the bilinear form (call it $\beta$) implies that for each $f$ in $\Hom{C}(\Sigma T_1^X,Y)$ and each $g$ in $\Hom{C}(Y, \Sigma^2 X)$, $\beta(\eta^*f, g) = \beta(f, \Sigma^2 \eta_* g)$.

As a consequence, there is an induced non-degenerate bilinear form
    \begin{displaymath}
	    \Ima \eta^* \times \Ima \Sigma^2 \eta_* \longrightarrow k.
    \end{displaymath}
    
Since $\Ima \eta^*$ is isomorphic to $(\Sigma T)(X,Y)$ and $\Ima \Sigma^2 \eta_*$ is in turn isomorphic to $\Hom{\cC}(Y, \Sigma^2 X)/(\Sigma^2 T)$, we get the desired result.

}

\subsection{Cluster character : definition}\label{sect::deficlusterchar}

In \cite{Palu08}, Y. Palu defined the notion of \emph{cluster character} for a $\Hom{}$--finite 2-Calabi--Yau triangulated category with a cluster-tilting object.  In our context, the category $\cC$ is not $\Hom{}$--finite nor 2-Calabi--Yau, and the object $T$ is only assumed to be rigid.  However, the definition can be adapted to this situation as follows.

\begin{definition}\label{defi::d}
Let $\cC$ be a triangulated category and $T$ be a rigid object as above.  The category $\cD$ is the full subcategory of $\pr{\cC}T \cap \pr{\cC}\Sigma T$ whose objects are those $X$ such that $FX$ is a finite-dimensional $B$--module.
\end{definition}

Under the hypotheses of this Section, the subcategory $\cD$ is Krull--Schmidt and stable under extensions.  Moreover, in the special case where $\cC = \cC_{Q,W}$ and $T = \Sigma^{-1}\Gamma$, the subcategory $\cD$ does not depend on the mutation class of $T$; that is, replacing $\Gamma$ by any $\mu_r\ldots\mu_1 \Gamma$ in the definition of $T$ yields the same subcategory $\cD$ (this is a consequence of the nearly Morita equivalence of \cite[Corollary 4.6]{KY09} and of Corollary \ref{coro::pr}).

The subcategory $\cD$ allows us to extend the notion of cluster characters.

\begin{definition}\label{defi::char}
Let $\cC$ be a triangulated category and $T$ be a rigid object as above.  Let $\cD$ be as Definition \ref{defi::d}.

A \emph{cluster character} on $\cC$ (with respect to $T$) with values in a commutative ring $A$ is a map
 \begin{displaymath}
  \chi : obj(\cD) \longrightarrow A	
 \end{displaymath}
satisfying the following conditions :
\begin{itemize}
	\item if $X$ and $Y$ are two isomorphic objects in $\cD$, then we have $\chi(X) = \chi(Y)$;
	\item for all objects $X$ and $Y$ of $\cD$, $\chi(X\oplus Y) = \chi(X)\chi(Y)$;
	\item (multiplication formula) for all objects $X$ and $Y$ of $\cD$ that are such that $\dim \Ext{1}{\cC}(X,Y) = 1$, the equality
	 \begin{displaymath}
     \chi(X)\chi(Y) = \chi(E) + \chi(E')
   \end{displaymath}
holds, where $X\longrightarrow E \longrightarrow Y \longrightarrow \Sigma X$ and $Y\longrightarrow E' \longrightarrow X \longrightarrow \Sigma Y$ are non split triangles.
\end{itemize}
\end{definition}

Note that $\Ext{1}{\cC}(Y,X)$ is one-dimensional, thanks to the non-degenerate bilinear form.  Also note that $E$ and $E'$ are in $\cD$, thanks to Lemma \ref{lemm::index2}, so $\chi(E)$ and $\chi(E')$ are defined.

\begin{remark}
If the category $\cC$ happens to be $\Hom{}$--finite and 2-Calabi--Yau, and if $T$ is a cluster-tilting object, then this definition is equivalent to the one given in \cite{Palu08}.  Indeed, in that case, it was shown in \cite{KR07}, Proposition 2.1, that $\pr{\cC}T = \pr{\cC}\Sigma T = \cC$, so $\cD = \cC$.
\end{remark}

Let $\cD$ be as in Definition \ref{defi::char} and $\iota$ be as in Notation \ref{nota::iota}.  Define the map 

\begin{displaymath}
	X'_{?} : obj(\cD) \longrightarrow \bQ(x_1, x_2, \ldots, x_n)
\end{displaymath}
as follows : for any object $X$ of $\cD$, put
\begin{displaymath}
	X'_{X} = x^{\footnotesize\ind{T}\Sigma^{-1}X} \sum_{e} \chi\Big(\Gr{e}(FX)\Big) x^{-\iota(e)}.
\end{displaymath}

Here, $\chi$ is the Euler--Poincar\'e characteristic.

Compare this definition to that of \cite{Palu08} and \cite{FK09}.

\begin{theorem}\label{theo::char}
The map $X'_{?}$ defined above is a cluster character on $\cC$ with respect to $T$.
\end{theorem}

It is readily seen that the first two conditions of Definition \ref{defi::char} are satisfied by $X'_{?}$.  We thus need to show that the multiplication formula holds in order to prove Theorem \ref{theo::char}.

\subsection{Dichotomy}\label{sect::dicho}

This subsection mimics Section 4 of \cite{Palu08}. Our aim here is to prove the following \emph{dichotomy} phenomenon.

Let $X$ and $Y$ be objects of $\cD$ such that $\dim \Ext{1}{\cC}(X, Y) = 1$. This implies that $\dim \Ext{1}{\cC}(Y, X) = 1$.  Let
\begin{displaymath}
	\xymatrix{X \ar[r]^{i} & E \ar[r]^{p} & Y \ar[r]^{\varepsilon} & \Sigma X \\
	          Y \ar[r]^{i'} & E' \ar[r]^{p'} & X \ar[r]^{\varepsilon'} & \Sigma Y}
\end{displaymath}
 
be non-split triangles.  Recall that Lemma \ref{lemm::index2} implies that $E$ and $E'$ are in $\cD$.

Let $U$ and $V$ be submodules of $FX$ and $FY$, respectively.  Define
\begin{displaymath}
	G_{U,V} = \Big\{ W \in \bigcup_{e}\Gr{e}(FE) \ \Big| \ (Fi)^{-1}(W) = U, \ Fp(W) = V \Big\} \textrm{ and}
\end{displaymath}
\begin{displaymath}
	G'_{U,V} = \Big\{ W \in \bigcup_{e}\Gr{e}(FE') \ \Big| \ (Fi')^{-1}(W) = V, \ Fp'(W) = U \Big\}.
\end{displaymath}

\begin{proposition}[Dichotomy]\label{prop::dicho}
Let $U$ and $V$ be as above.  Then exactly one of $G_{U,V}$ and $G'_{U,V}$ is non-empty.
\end{proposition}

In order to prove this Proposition, a few lemmata are needed.

Using Lemma \ref{lemm::modules}, lift the inclusions $U \subseteq FX$ and $V \subseteq FY$ to morphisms $i_U:\overline{U}\longrightarrow X$ and $i_V:\overline{V}\longrightarrow Y$, where $\overline{U}$ and $\overline{V}$ lie in $\pr{\cC}T \cap \pr{\cC}\Sigma T$.  Keep these notations for the rest of this Section and for the next.

The first lemma is about finiteness.

\begin{lemma}\label{lemm::finext}
Let $X$ and $\overline{U}$ be as above. Let $M$ be an object of $\cC$ such that $FM$ and $\Hom{\cC}(X, \Sigma M)$ are finite-dimensional.  Then $\Hom{\cC}(\overline{U}, \Sigma M)$ is also finite-dimensional.
\end{lemma}
\demo{Embed $i_U$ in a triangle $\xymatrix{\Sigma^{-1} X \ar[r]^{\phantom{xx}\alpha} & Z \ar[r]^{\beta} & \overline{U} \ar[r]^{i_U} & X}$.  From this triangle, one gets the exact sequence
\begin{displaymath}
	\xymatrix{\Hom{\cC}(X, \Sigma M)\ar[r]^{i_U^*} & \Hom{\cC}(\overline{U}, \Sigma M)\ar[r]^{\beta^*} & \Hom{\cC}(Z, \Sigma M).}
\end{displaymath}
The image of $\beta^*$ is isomorphic to $\Hom{\cC}(\overline{U}, \Sigma M)/\Ima i_U^*$.  Since $\Hom{\cC}(X, \Sigma M)$ is finite-dimensional by hypothesis, it suffices to show that $\Ima \beta^*$ is finite-dimensional to prove the Lemma.

Since $\Sigma^{-1} X$ and $\overline{U}$ are in $\pr{\cC}T$ and $\Coker F\beta$ is finite-dimensional, Lemma \ref{lemm::index2} can be applied to get that $Z$ is in $\pr{\cC}T$.  Let $\xymatrix{T_1^Z\ar[r] & T_0^Z \ar[r]& Z \ar[r] & \Sigma T_1^Z}$ be a triangle, with $T_0^Z$ and $T_1^Z$ in $\pr{\cC} T$.

Now, $Fi_U$ is a monomorphism, so $F\beta = 0$, and Lemma \ref{lemm::modules} implies that $\beta$ lies in $(\Sigma T)$.  Therefore $\Ima \beta^*$ is contained in $(\Sigma T)(Z, \Sigma M)$.  It is thus sufficient to show that the latter is finite-dimensional.

We have an exact sequence
\begin{displaymath}
	\xymatrix{\Hom{\cC}(\Sigma T_1^Z, \Sigma M)\ar[r]^{\gamma} & \Hom{\cC}(Z, \Sigma M)\ar[r] & \Hom{\cC}(T_0^Z, \Sigma M).}
\end{displaymath}
Since $T$ is rigid, we have that $(\Sigma T)(Z, \Sigma M) = \Ima \gamma$.  It is thus finite-dimensional. Indeed, $FM = \Hom{\cC}(T, M)$ is finite-dimensional, and this implies that the same property holds for $\Hom{\cC}(\Sigma T_1^Z, \Sigma M)$.  This finishes the proof.
}

The second lemma is a characterization in $\cC$ of the non-emptiness of $G_{U,V}$.  It is essentially \cite[Lemma 14]{Palu08}, where the proof differs in that not every object lies in $\pr{\cC}T$.

\begin{lemma}\label{lemm::nonempty}
With the above notations, the following are equivalent:
\begin{enumerate}
  \item $G_{U,V}$ is non-empty;
  \item there exist morphisms $e:\Sigma^{-1}\overline{V} \longrightarrow \overline{U}$ and $f:\Sigma^{-1}Y \longrightarrow \overline{U}$ such that
    \begin{enumerate}
      \item $(\Sigma^{-1}\varepsilon)(\Sigma^{-1}i_V) = i_Ue$
      \item $e \in (T)$
      \item $i_Uf - \Sigma^{-1}\varepsilon \in (\Sigma T)$;
    \end{enumerate}
  \item condition (2) where, moreover, $e = f\Sigma^{-1}i_V$.
\end{enumerate}
\end{lemma}
\demo{Let us first prove that (2) implies (1).  The commutative square given by (a) gives a morphism of triangles
\begin{displaymath}
	\xymatrix{ \overline{U}\ar[r]\ar[d]^{i_U} & \overline{W}\ar[r]\ar[d]^{\phi} & \overline{V}\ar[r]^{\Sigma e}\ar[d]^{i_V} & \Sigma\overline{U}\ar[d]^{\Sigma i_U} \\
	           X\ar[r]^{i} & E\ar[r]^{p} & Y\ar[r]^{\varepsilon} & \Sigma X.
	}
\end{displaymath}

Applying the functor $F$, we get a commutative diagram in $\MOD B$ :

\begin{displaymath}
	\xymatrix{ U\ar[r]\ar[d]^{Fi_U} & F\overline{W}\ar[r]\ar[d]^{F\phi} & V\ar[r]^{\Sigma e}\ar[d]^{Fi_V} & 0 \\
	           FX\ar[r]^{Fi} & FE\ar[r]^{Fp} & FY. & 
	}
\end{displaymath}

An easy diagram chasing shows that the image of $F\phi$ is in $G_{U,V}$, using the morphism $f$.

Let us now prove that (1) implies (2).  Let $W$ be in $G_{U,V}$.  Notice that $U$ contains $\Ker Fi = \Ima F\Sigma^{-1}\varepsilon$, so $F\Sigma^{-1}\varepsilon$ factors through $U$.  Since $\Sigma^{-1}Y \in \pr{\cC}T$, Lemma \ref{lemm::modules} allows us to find a lift $f:\Sigma^{-1}Y \longrightarrow \overline{U}$ of this factorization such that $i_Uf - \Sigma^{-1}\varepsilon \in  (\Sigma T)$.

Let us define $e$.  Since $\overline{V}\in\pr{\cC}T$, there exists a triangle

\begin{displaymath} 
\xymatrix{ T_1^V \ar[r] & T_0^V \ar[r]& \overline{V}\ar[r] & \Sigma T_1^V}. 
\end{displaymath}

Since $FT_0^V$ is projective, and since $\xymatrix{W\ar[r]^{Fp} & V}$ is surjective,  $\xymatrix{FT_0^V\ar[r] & V}$ factors through $Fp$.  Composing with the inclusion of $W$ in $FE$, we get a commutative square
\begin{displaymath}
	\xymatrix{ FT_0^V \ar[r]\ar[d] & V\ar[d] \\
	           FE     \ar[r]       & FY
	}
\end{displaymath}
which lifts to a morphism of triangles (thanks to Lemma \ref{lemm::modules})
\begin{displaymath}
	\xymatrix{ \Sigma^{-1}\overline{V} \ar[r]\ar[d] & T_1^V\ar[r]\ar[d] & T_0^V\ar[r]\ar[d] & \overline{V}\ar[d]  \\
	           \Sigma^{-1}Y\ar[r] & X \ar[r] & E \ar[r] & Y.
	}
\end{displaymath}
Now $\xymatrix{FT_0^V\ar[r] & FE}$ factors through $W$, and since $U = (Fi)^{-1}(W)$, then the image of $\xymatrix{FT_1^V \ar[r] & FX}$ is contained in $U$.  Thus $\xymatrix{T_1^V \ar[r] & X}$ factors through $\overline{U}$, and we take $e$ to be the composition $\xymatrix{\Sigma^{-1}V \ar[r] & T_1^V \ar[r] & U}$.  By construction, conditions (a) and (b) are satisfied.

Obviously, (3) implies (2).  Let us show that (2) implies (3).  

First, since $(\Sigma^{-1}\varepsilon)(\Sigma^{-1}i_V) = i_Ue$ and $i_Uf\Sigma^{-1}i_V - (\Sigma^{-1}\varepsilon)(\Sigma^{-1}i_V) \in (\Sigma T)$, we get that $i_U(f\Sigma^{-1}i_V - e) \in (\Sigma T)$.  Since $Fi_U$ is a monomorphism, and since $\Sigma^{-1}V \in \pr{\cC}T$, we get that $h:=f\Sigma^{-1}i_V - e \in (\Sigma T)$.

Embed the morphism $\Sigma^{-1}i_V$ into a triangle 

\begin{displaymath}
\xymatrix{\Sigma^{-1}\overline{V}\ar[r]^{\Sigma^{-1}i_V} & \Sigma^{-1}Y\ar[r] & C\ar[r] & \overline{V}}.
\end{displaymath}
  Using Lemma \ref{lemm::index2}, we see that $C$ lies in $\pr{\cC}T$, and since $Fi_V$ is a monomorphism, this implies that $\xymatrix{C \ar[r] & \overline{V}}$ lies in $(\Sigma T)$, by Lemma \ref{lemm::modules}.

Now, $h\in (\Sigma T)$ and $\xymatrix{\Sigma^{-1}C \ar[r] & \Sigma^{-1}\overline{V}} \in (T)$, so their composition vanishes.  Therefore there exists a morphism $\xymatrix{\Sigma^{-1}Y \ar[r]^{\ell} & \overline{U}}$ such that $\ell\Sigma^{-1}i_V = h$.

Since $\overline{V}$ is in $\pr{\cC}\Sigma T$, there is a triangle $\xymatrix{T_V^1\ar[r] & T_V^0\ar[r] & \Sigma^{-1}\overline{V}\ar[r] & \Sigma T_V^1}$.  Now, since $\xymatrix{\Sigma^{-1}C \ar[r] & \Sigma^{-1}\overline{V}} \in (T)$ and $\ell\Sigma^{-1}i_V \in (\Sigma T)$, we have morphisms of triangles

\begin{displaymath}
	\xymatrix{ \Sigma^{-1}C \ar[r]\ar[d] & \Sigma^{-1}\overline{V}\ar[r]^{\Sigma^{-1}i_V}\ar@{=}[d] & \Sigma^{-1}Y\ar[r]^{c}\ar[d]^{v} & C\ar[d]  \\
	           T_V^0\ar[r]^{u}\ar[d] & \Sigma^{-1}\overline{V}\ar[r]\ar[d]^{\Sigma^{-1}i_V} & \Sigma T_V^1\ar[r]\ar[d]^{w} & \Sigma T_V^0\ar[d] \\
	           \Sigma^{-1}C'\ar[r] & \Sigma^{-1}Y\ar[r]^{\ell} & \overline{U}\ar[r] & C'.
	}
\end{displaymath}

Since the composition $(\ell - wv)\Sigma^{-1}i_V$ vanishes, there exists a morphism $\ell'$ from $C$ to $\overline{U}$ such that $\ell'c = \ell - wv$.

Put $f_0 = f - wv$.  Then
\begin{displaymath}
	f_0\Sigma^{-1}i_V = f\Sigma^{-1}i_V - \ell\Sigma^{-1}i_V + \ell'c\Sigma^{-1}i_V = f\Sigma^{-1}i_V - h + 0 = e.  
\end{displaymath}
Moreover,  $i_Uf_0 - \Sigma^{-1}\varepsilon = (i_Uf - \Sigma^{-1}\varepsilon) - wv \in (\Sigma T)$.  This finishes the proof.
}

\demo{(of Proposition \ref{prop::dicho}) The proof is similar to that of Proposition 15 of \cite{Palu08}.  Consider the linear map
\begin{displaymath}
  \begin{array}{ccc}
  \alpha: \! (\Sigma^{-1}Y, \overline{U}) \oplus (\Sigma^{-1}Y, X)  & \rightarrow & \! (\Sigma^{-1}\overline{V}, X) \oplus (\Sigma^{-1}\overline{V}, \overline{U})/(T) \oplus (\Sigma^{-1}Y, X)/(\Sigma T)\\
  (x,y) & \mapsto & ( y(\Sigma^{-1}i_V) - i_Ux(\Sigma^{-1}i_V), \ x\Sigma^{-1}i_V, i_Ux - y), 
  \end{array} 
\end{displaymath}
where we write $(X,Y)$ instead of $\Hom{\cC}(X,Y)$.  Then $f \in \Hom{\cC}(\Sigma^{-1}Y, \overline{U})$ satisfies condition (3) of Lemma \ref{lemm::nonempty} if, and only if, $(f, \Sigma^{-1}\varepsilon)$ is in $\Ker \alpha$.  Since $\Hom{\cC}(Y, \Sigma^{-1}X)$ is one-dimensional, the existence of such an $f$ is equivalent to the statement that the map
\begin{displaymath}
	\xymatrix{\beta:\Ker \alpha \ar@{^{(}->}[r] & (\Sigma^{-1}Y, \overline{U}) \oplus (Y, \Sigma X) \ar@{->>}[r] & (Y, \Sigma X) }
\end{displaymath} 
does not vanish, where the second map is the canonical projection.

Now, the emptiness of $G_{U,V}$ is equivalent to the vanishing of $\beta$, which is equivalent to the vanishing of its dual
\begin{displaymath}
	\xymatrix{D\beta: D(Y, \Sigma X) \ar@{^{(}->}[r] & D(\Sigma^{-1}Y, \overline{U}) \oplus D(Y, \Sigma X) \ar@{->>}[r] &  \Coker D\alpha, }
\end{displaymath}
which is in turn equivalent to the fact that any element of the space $D(\Sigma^{-1}Y, \overline{U}) \oplus D(Y, \Sigma X)$ of the form $(0,z)$ lies in $\Ima D\alpha$.

Using Lemma \ref{lemm::finext} and the non-degenerate bilinear form, we see that all five spaces involved in the definition of $\alpha$ are finite-dimensional.  Therefore, Lemma \ref{lemm::bil} yields the following commutative diagram, where the vertical morphisms are componentwise isomorphisms :
\begin{displaymath}
  \xymatrix{
    D(\Sigma^{-1}\overline{V}, X) \oplus D(\Sigma^{-1}\overline{V}, \overline{U})/(T) \oplus D(\Sigma^{-1}Y, X)/(\Sigma T)\ar[r]^{\phantom{xxxxxxxxxxxx}D\alpha}   &  D(\Sigma^{-1}Y, \overline{U}) \oplus D(\Sigma^{-1}Y, X)\\
    (X, \Sigma\overline{V}) \oplus (\Sigma T)(\overline{U}, \Sigma\overline{V}) \oplus (\Sigma^2 T)(X, \Sigma Y) \ar[r]^{\phantom{xxxxxxxxxxxx}\alpha'}\ar[u]   &  (\overline{U}, \Sigma Y) \oplus (X, \Sigma Y).\ar[u]   
  } 
\end{displaymath}

Here, the action of $\alpha'$ is given by $\alpha'(x,y,z) = ((\Sigma i_V)xi_U + (\Sigma i_V)y + zi_U, -(\Sigma i_V)x -z)$.

Now, any element of the form $(0,w)$ is in $\Ima \alpha'$ if, and only if, $(0, \varepsilon')$ is in $\Ima \alpha'$, which in turn is equivalent to the fact that condition (2) of Lemma \ref{lemm::nonempty} is satisfied, and thus to the fact that $G'_{U,V}$ is non-empty.  This finishes the proof.
}

\subsection{Multiplication formula}\label{sect::multiplication}
The main goal of this section is to prove the following Proposition, using the results of the previous sections.

\begin{proposition}[Multiplication formula]\label{prop::mult}
The map $X'_{?}$ satisfies the multiplication formula given in Definition \ref{defi::char}.
\end{proposition}

In order to prove this result, some notation is in order.  Let $X$ and $Y$ be objects of $\cD$ such that $\Hom{\cC}(X, \Sigma Y)$ is one-dimensional.  Let

\begin{displaymath}
	\xymatrix{ X\ar[r]^{i} & E\ar[r]^{p} & Y\ar[r]^{\varepsilon} & \Sigma X \textrm{ and} \\
	           Y\ar[r]^{i'} & E'\ar[r]^{p'} & X\ar[r]^{\varepsilon'} & \Sigma X
	}
\end{displaymath}

be non-split triangles in $\cC$.  For any submodules $U$ of $FX$ and $V$ of $FY$, define $G_{U,V}$ and $G'_{U,V}$ as in Section \ref{sect::dicho}.

For any dimension vectors $e$, $f$ and $g$, define the following varieties :

$$G_{e,f} = \bigcup_{\scriptsize{\begin{array}{c} 
	                    \dimv U = e \\ 
	                    \dimv V = f
	                   \end{array}}}  G_{U,V}$$ 
 	                   
	$$G'_{e,f} = \bigcup_{\scriptsize{\begin{array}{c} 
	                    \dimv U = e \\ 
	                    \dimv V = f
	                   \end{array}}}  G'_{U,V}$$ 
	                   
	$$G_{e,f}^{g} = G_{e,f} \cap \Gr{g}(FE)$$ 
  	
	$$G_{e,f}^{'g} = G'_{e,f} \cap \Gr{g}(FE').$$

We first need an equality on Euler characteristics.

\begin{lemma}\label{lemm::Euler}
  With the above notation, we have that 
  \begin{displaymath}
   \chi(\Gr{e}(FX))\chi(\Gr{f}(FY)) = \sum_{g}\Big( \chi(G_{e,f}^{g}) + \chi(G_{e,f}^{'g}) \Big).
  \end{displaymath} 
\end{lemma}
\demo{The Lemma is a consequence of the following successive equalities:
\begin{eqnarray*}
	\chi(\Gr{e}(FX))\chi(\Gr{f}(FY)) & = & \chi(\Gr{e}(FX) \times \Gr{f}(FY)) \\
	                                 & = & \chi(G_{e,f} \sqcup G'_{e,f}) \\
	                                 & = & \chi(G_{e,f}) + \chi(G'_{e,f}) \\
	                                 & = & \sum_{g} \Big( \chi(G_{e,f}^{g}) + \chi(G_{e,f}^{'g})   \Big).	
\end{eqnarray*} 

The only equality requiring explanation is the second one.  Consider the map
\begin{eqnarray*}
  G_{e,f} \sqcup G'_{e,f} & \longrightarrow & \Gr{e}(FX) \times \Gr{f}(FY) \\
  W & \longmapsto & \left\{ \begin{array}{ll} 
                               \Big( (Fi)^{-1}(W), (Fp)(W) \Big) & \textrm{if } W \in G_{e,f}	 \\ 
                               \Big( (Fi')^{-1}(W), (Fp')(W) \Big) & \textrm{if } W \in G'_{e,f}. 
                            \end{array} \right.
\end{eqnarray*}

As a consequence of Proposition \ref{prop::dicho}, the map is surjective, and as shown in \cite{CC06}, its fibers are affine spaces.  The Euler characteristic of all its fibers is thus $1$, and we have the desired equality. 
}

Secondly, we need an interpretation of the dimension vectors $e$, $f$ and $g$.

\begin{lemma}\label{lemm::dimv}
If $G_{e,f}^{g}$ is not empty, then
 \begin{displaymath}
	\dimv (\Coker F\Sigma^{-1}p) = e + f - g.
 \end{displaymath}
\end{lemma}
\demo{ We have the following commuting diagram in $\MOD B$, where the rows are exact sequences:
 \begin{displaymath}
  \xymatrix{ 0\ar[r] & K\ar[r] & U\ar[r]\ar[d]^{Fi_U} & W\ar[r]\ar[d]^{Fi_W} & V\ar[r]\ar[d]^{Fi_V} & 0 \\
             & & FX\ar[r]^{Fi} & FE\ar[r]^{Fp} & FY\ar[r]^{F\varepsilon} & F\Sigma X.
  }	
 \end{displaymath}
 
In this diagram, $\dimv U = e$, $\dimv W = g$ and $\dimv V = f$.  The existence of such a diagram is guaranteed by the non-emptiness of $G_{e,f}^{g}$.

Now $\Coker F\Sigma^{-1}p$ is isomorphic to $\Ker Fi$, which is in turn isomorphic to $\Ker (Fi \circ Fi_U)$, since $U = (Fi)^{-1}(W)$.  This last kernel is isomorphic to $K$.  Hence the equality $\dimv (\Coker F\Sigma^{-1}p) = \dimv K$ holds.

Finally, the upper exact sequence gives the equality $\dimv K + \dimv W = \dimv U + \dimv V$.  By rearranging and substituting terms, we get the desired equality.
}

Everything is now in place to prove the multiplication formula.

\demo{ (of Proposition \ref{prop::mult}) The result is a consequence of the following successive equalities (explanations follow).
  \begin{eqnarray*}
	  X'_X X'_Y & = & x^{\footnotesize\ind{T} \Sigma^{-1}X + \ind{T} \Sigma^{-1}Y} \sum_{e,f} \chi\Big( \Gr{e}(FX) \Big) \chi\Big( \Gr{f}(FY) \Big) x^{-\iota(e+f)} \\
	  & = & x^{\footnotesize\ind{T} \Sigma^{-1}X + \ind{T} \Sigma^{-1}Y} \sum_{e,f,g} \Big( \chi(G_{e,f}^{g}) + \chi(G_{e,f}^{'g}) \Big) x^{-\iota(e+f)}\\
	  & = &  x^{\footnotesize\ind{T} \Sigma^{-1}X + \ind{T} \Sigma^{-1}Y -  \iota(\Coker F\Sigma^{-1}p) - \iota(g)} \sum_{e,f,g} \chi(G_{e,f}^{g}) + \nonumber \\
	  &   & {} + \ x^{\footnotesize\ind{T} \Sigma^{-1}X + \ind{T} \Sigma^{-1}Y -  \iota(\Coker F\Sigma^{-1}p') - \iota(g)} \sum_{e,f,g} \chi(G_{e,f}^{'g}) \\
	  & = & x^{\footnotesize\ind{T} \Sigma^{-1}E} \sum_{g} \chi\Big(\Gr{g}(FE)\Big)x^{-\iota (g)} + \nonumber \\
	  &   & {} + \ x^{\footnotesize\ind{T} \Sigma^{-1}E'} \sum_{g} \chi\Big(\Gr{g}(FE')\Big)x^{-\iota (g)} \\
	  & = & X'_E + X'_{E'}.
  \end{eqnarray*}	

The first equality is just the definition of $X'_?$.  The second one is a consequence of Lemma \ref{lemm::Euler}, and the third one of Lemma \ref{lemm::dimv}.  The fourth follows from Lemma \ref{lemm::index2}.  The fifth equality is obtained by definition of $G_{e,f}^{g}$ and $G_{e,f}^{'g}$.  The final equality is, again, just the definition of $X'_?$.
}
\section{Application to cluster algebras}\label{sect::application}

In this section, we apply the cluster character developped in Section \ref{sect::clusterchar} to any skew-symmetric cluster algebra, taking $T = \Sigma^{-1} \Gamma$.  

An object $X$ of a triangulated category is \emph{rigid} if $\Hom{}(X, \Sigma X)$ vanishes.  Let $\cC_{Q,W}$ be the cluster category of a quiver with potential $(Q,W)$.  A \emph{reachable object} of $\cC_{Q,W}$ is a direct factor of a direct sum of copies $\mu_{i_r}\ldots\mu_{i_1}(\Gamma)$ for some admissible sequence of vertices $(i_1, \ldots, i_r)$ of $Q$.  Notice that each reachable object is rigid.

\begin{theorem}\label{theo::main}
Let $(Q,W)$ be a quiver with potential.  Then the cluster character $X'_?$ defined in \ref{defi::char} gives a surjection from the set of isomorphism classes of indecomposable reachable objects of $\cC_{Q,W}$ to the set of clusters variables of the cluster algebra associated with $Q$ obtained by mutating the initial seed at admissible sequences of vertices. 

More precisely, $X'_?$ sends the indecomposable summands of $\mu_{i_r}\ldots\mu_{i_1}(\Gamma)$ to the elements of the mutated cluster $\mu_{i_r}\ldots\mu_{i_1}(x_1, \ldots, x_r)$, where $(x_1, \ldots, x_r)$ is the initial cluster.
\end{theorem}
\demo{ Let $(i_1, \ldots, i_r)$ be an admissible sequence of vertices.  It is easily seen that $X'_{\Gamma_i} = x_i$ for all vertices $i$.
It is a consequence of \cite[Corollary 4.6]{KY09} that $\Hom{\cC}(\Sigma^{-1}\Gamma, \mu_{i_r}\ldots\mu_{i_1}(\Gamma))$ is finite dimensional.  Moreover, $\mu_{i_r}\ldots\mu_{i_1}(\Gamma)$ is obviously in $\pr{\cC}(\mu_{i_r}\ldots\mu_{i_1}(\Gamma)) \cap \pr{\cC}(\Sigma^{-1}\mu_{i_r}\ldots\mu_{i_1}(\Gamma))$, which is equal to $\pr{\cC}\Gamma \cap \pr{\cC}\Sigma^{-1}\Gamma$ by Corollary \ref{coro::pr}.  Therefore $\mu_{i_r}\ldots\mu_{i_1}(\Gamma)$ is in the subcategory $\cD$ of Definition \ref{defi::d}, and we can apply $X'_?$ to $\mu_{i_r}\ldots\mu_{i_1}(\Gamma)$.

We prove the result by induction on $r$.  

First notice that $\mu_{i_r}\ldots\mu_{i_1}(\Gamma)_i = \mu_{i_{r-1}}\ldots\mu_{i_1}(\Gamma)_i$ if $i\neq i_r$.  Now, for $i = i_r$, using the triangle equivalence $\cC_{\mu_{i_r}\ldots\mu_{i_1}(Q,W)} \longrightarrow \cC_{Q,W}$, we get triangles
\begin{displaymath}
	\mu_{i_{r-1}}\ldots\mu_{i_1}(\Gamma)_{i_r} \longrightarrow \bigoplus_{i_r \rightarrow j}\mu_{i_{r-1}}\ldots\mu_{i_1}(\Gamma)_j \longrightarrow \mu_{i_{r}}\ldots\mu_{i_1}(\Gamma)_{i_r} \longrightarrow \ldots \end{displaymath}
\begin{displaymath}	
	\mu_{i_{r}}\ldots\mu_{i_1}(\Gamma)_{i_r} \longrightarrow \bigoplus_{j \rightarrow i_r}\mu_{i_{r-1}}\ldots\mu_{i_1}(\Gamma)_j \longrightarrow \mu_{i_{r-1}}\ldots\mu_{i_1}(\Gamma)_{i_r} \longrightarrow \ldots
\end{displaymath}
to which we can apply the multiplication formula of Proposition \ref{prop::mult}.  In this way, we obtain the mutation of variables in the cluster algebra.  This proves the result.
}

{\footnotesize
{}  
}

\end{document}